\documentclass[a4paper,twoside]{amsart}
\usepackage{amsmath}
\usepackage{amsfonts}
\usepackage{amssymb}
\usepackage{amsthm}
\usepackage{newlfont}
\usepackage{graphicx}
\usepackage{amscd}
\usepackage{color}
\usepackage{rotating}
\usepackage{genyoungtabtikz}
\def\shuffle{\begin{picture}(16,5)(-2,0) 
\put(0,0){\line(0,1){5}}
\put(5,0){\line(0,1){5}}
\put(10,0){\line(0,1){5}}
\put(0,0){\line(1,0){10}}
\end{picture}}
\def\Qshuffle{\stackrel{Q}{\shuffle}}
\def\Fshuffle{\stackrel{F}{\shuffle}}
\def\Tshuffle{\stackrel{T}{\shuffle}}

\def\sshuffle{\begin{picture}(8,3)(-1,0) 
\put(0,0){\line(0,1){3}}
\put(3,0){\line(0,1){3}}
\put(6,0){\line(0,1){3}}
\put(0,0){\line(1,0){6}}
\end{picture}}
\def\sFshuffle{\stackrel{F}{\sshuffle}}

\def\ttree{\begin{picture}(8,3) \put(3,2){\circle*{5}} \end{picture}}
\theoremstyle{plain}
\newtheorem{Th}{Theorem}[section]
\newtheorem{Cor}[Th]{Corollary}

\newtheorem{Prop}[Th]{Proposition}
\theoremstyle{definition}

\newtheorem{Ex}{Example}[section]
\theoremstyle{remark}
\newtheorem*{Rem}{Remark}%[section]
\numberwithin{equation}{section}
\begin{document}

\title[Generalized quasi-symmetric functions]
{Hopf algebra structure of generalized quasi-symmetric functions in partially commutative variables}

\author{Adam Doliwa}

\address{Adam Doliwa, Faculty of Mathematics and Computer Science, University of Warmia and Mazury in Olsztyn,
ul.~S{\l}oneczna~54, 10-710~Olsztyn, Poland}

\email{doliwa@matman.uwm.edu.pl}
\urladdr{http://wmii.uwm.edu.pl/~doliwa/}

\date{}
\begin{abstract}
We introduce a coloured generalization  $\mathrm{NSym}_A$ of the Hopf algebra of non-commutative symmetric functions  described as a subalgebra of the of rooted ordered coloured trees Hopf algebra. Its natural basis can be identified with the set of sentences over alphabet $A$ (the set of colours). We present also its graded dual algebra $\mathrm{QSym}_A$ of coloured quasi-symmetric functions together with its realization in terms of power series in partially commutative variables.  We provide formulas expressing multiplication, comultiplication and the antipode for these Hopf algebras in various bases --- the corresponding generalizations of the complete homogeneous, elementary, ribbon Schur and power sum bases of $\mathrm{NSym}$, and the monomial and fundamental bases of $\mathrm{QSym}$. We study also certain distinguished series of trees in the setting of restricted duals to Hopf algebras.
\end{abstract}
\maketitle

\section{Introduction} 
Theory of Hopf algebras forms a modern basis for understanding symmetries of solvable models in quantum and statistical theoretical physics~\cite{Majid,Klimyk-Schmudgen,Chari-Pressley}.
Application of Hopf algebras~\cite{Abe,Sweedler} to combinatorics can be traced back to Rota~\cite{Joni-Rota}, see also ~\cite{Schmitt,Grinberg-Reiner} for more recent reviews of the subject, which has expanded since then.  Combinatorial aspects of the Bethe ansatz and of the quantum inverse scattering method \cite{QISM} were studied, for example, in works by Fomin, Kirillov and Reshetikhin~\cite{Fomin-Kirillov,Kirillov-Reshetikhin}. For more about mutual interactions between the theory of integrable systems and combinatorics, see recent reviews \cite{Deift,tau-generating,diFrancesco,Zinn-Justin}.

Hopf algebras of rooted trees appeared in the analysis of Runge--Kutta methods by Butcher~\cite{Butcher} and D\"{u}r~\cite{Dur}, and in works by Grossman and Larson~\cite{Grossman-Larson} in the context of symbolic computation. 
More recently they were used by Connes and Kreimer~\cite{Connes-Kreimer} to describe renormalization procedure of quantum field theory, see also 
\cite{Broadhurst-Kreimer,Brouder-Frabetti}. The non-commutative Hopf algebra of trees and forests, generalizing that of Connes and Kreimer, was considered by Foissy~\cite{Foissy-I}, and independently by Holtkamp~\cite{Holtkamp}. 

The theory of symmetric functions~\cite{Stanley,Macdonald} is by now well established subject with numerous applications in algebraic topology, combinatorics, representation theory, integrable systems and geometry. Quasi-symmetric functions, introduced by Gessel~\cite{Gessel-q-sym} (see also an earlier relevant work of Stanley~\cite{Stanley-q-sym}), are extensions of symmetric functions that are becoming of comparable importance~\cite{Introduction-QSym,Baker-Richter,Mason-qsym}. As a graded Hopf algebra, the dual of the algebra of quasi-symmetric functions is the Hopf algebra of non-commutative symmetric functions 
introduced by Gelfand, Krob, Lascoux, Leclerc, Retakh, and Thibon~\cite{Gelfand-NSym}. In works of Zhao \cite{Zhao} and Hoffman~\cite{Hoffman-nc} there was established isomorphism between the Hopf algebra of non-commutative symmetric functions and certain subalgebra of rooted ordered trees, see also~\cite{EbrahimiFard-Guo-Kreimer} for application of such trees (called ladders) to study integrable aspect of the renormalization.

Our work arose from the search for generalization of the relationship between theory of symmetric functions, combinatorics and the integrable systems on the non-commutative level. Already the standard description of the Kadomtsev--Petviashvili (KP) hierarchy of integrable partial differential equations in terms of free fermions by the Kyoto School~\cite{Miwa-Jimbo-Date} involves large part of the theory of symmetric functions~\cite{Sato}, see also a generalization \cite{Dimakis-Muller-Hoisen} in direction of quasi-symmetric functions. For example, the Schur functions when expressed in suitably scaled power sum functions (times of the KP hierarchy) provide polynomial $\tau$-function solutions of the equations.
Non-commutative extensions of integrable systems are of growing interest in mathematical physics~\cite{Kupershmidt,BobSur-nc,Nimmo-NCKP,EtingofGelfandRetakh,DF-K,Konstantinou-Rizos,Dol-Des,Doliwa-YB,DoliwaNoumi,DoliwaKashaev}. 
In this paper we define and study properties of a coloured version of the Hopf algebra of non-commutative symmetric function and of its graded dual.  The idea to consider coloured versions of various algebras is not new, see for example~\cite{BaumannHohlweg,Hsiao,HsiaoPetersen,MantaciReutenauer,NovelliThibon,NovelliThibon2,Poirier}, where some generalizations of the Hopf algebras of non-commutative symmetric functions or quasi-symmetric functions have been discussed as well. The generalization presented in our work is to our best knowledge new, and in particular it extends some of the previous concepts, see  the last remark of Section~\ref{sec:sent-A}.

Let us present the structure of the paper, where we present step by step our generalization of the basic structural elements of the theory of quasi-symmetric and non-commutative symmetric functions showing similarities and differences with the original theory. In introductory Section~\ref{sec:Hopf} we recall necessary elements of the theory of graded Hopf algebras. As basic example we take the free Hopf algebra over finite alphabet, the Hopf algebra of quasi-symmetric functions, and their duals --- the shuffle algebra and the algebra of non-commutative symmetric functions. We pay also special attention to the Hopf algebra of rooted ordered coloured trees, closely related to the algebraic renormalization procedure  of the quantum field theory. 

Then in Section~\ref{sec:n-NSyn-QSym} we study in detail the Hopf algebra of sentences (coloured compositions or tall trees), the partial order in the set of sentences, and we provide briefly an interpretation of the sentences as words of certain basic context-free language. Then we change slightly our point of view by presenting the sentences over alphabet $A$ as the analog of the complete homogeneous basis in our coloured non-commutative generalization $\mathrm{NSym}_A$ of the Hopf algebra of symmetric functions. Then we introduce the corresponding analog of the basis of elementary functions and discuss coloured version of the standard formulas describing mutual interrelation between the complete homogeneous and elementary functions. 

We devote Section~\ref{sec:CQSym} to description of the Hopf algebra $\mathrm{QSym}_A$ of coloured quasi-symmetric functions as graded dual of $\mathrm{NSym}_A$. We define first the basis of coloured monomial quasi-symmetric functions in the standard way as the dual basis to the coloured complete homogeneous functions. Then we construct its `polynomial' realization in terms of certain power series of bounded degree in partially commuting variables. Such a partial commutativity is completely new ingredient of our generalization of the theory of quasi-symmetric functions. We remark that partially commutative variables have been introduced to study combinatorial problems by Cartier and Foata in \cite{Cartier-Foata}. They also have found applications in algebra, theory of orthogonal polynomials, statistical physics and computer science; see review by Viennot~\cite{Viennot-heaps} written in terms of heaps of pieces. In theoretical computer science, as was proposed by Mazurkiewicz \cite{Mazurkiewicz}, they describe concurrent computations. We would like to stress that both algebras  $\mathrm{NSym}_A$ and  $\mathrm{QSym}_A$ are non-commutative and non-cocommutative for $|A|>1$.

Sections~\ref{sec:F-R} and~\ref{sec:restricted} are devoted to presentation of the coloured generalization of other rudimentary elements of the theory of symmetric functions. We define and study the fundamental basis of $\mathrm{QSym}_A$ and its dual basis in $\mathrm{NSym}_A$ of coloured non-commutative ribbon Schur functions. Finally we construct the coloured non-commutative version of the power sum symmetric functions. In Section~\ref{sec:restricted} we extended also the algebra of trees to its infinite series version working within the setting of the restricted dual Hopf algebras. This point of view was useful in studying the coloured non-commutative power sum functions, but certainly deserves deeper studies in the context of applications of combinatorial Hopf algebras to integrable systems and theoretical physics.

\subsection*{Acknowledgements} 
%I would like to thank to a Referee for careful reading the paper. 
The research was supported by National Science Centre, Poland, under grant 2015/19/B/ST2/03575 \emph{Discrete integrable systems -- theory and applications}.

\section{Hopf algebras of trees and quasi-symmetric functions} \label{sec:Hopf}
We assume that the Reader is familiar with the basic definitions and properties of Hopf algebras, as covered in \cite{Abe} or \cite{Sweedler}. All the results presented in this Section are known, but we recall them to provide necessary terminology and background to formulate new ones in the next Sections. In the paper all algebras are over a fixed field $\Bbbk$ of characteristic zero, although sometimes a commutative ring may be enough.

\subsection{Hopf algebras}
By $(\mathcal{H},\mu,\eta,\Delta,\epsilon)$ denote a bialgebra which is:
\begin{enumerate}
\item an associative algebra $(\mathcal{H},\mu,\eta)$ consisting of $\Bbbk$-linear
multiplication $\mu\colon \mathcal{H}\otimes \mathcal{H} \to \mathcal{H}$ and $\Bbbk$-linear unit map $\eta\colon \Bbbk \to \mathcal{H} $ satisfying properties described by the diagrams:
\begin{equation}
\begin{CD}
\mathcal{H}\otimes \mathcal{H}\otimes \mathcal{H} @ > \mathrm{id} \otimes \mu >>  \mathcal{H} \otimes \mathcal{H} \\
@V \mu \otimes \mathrm{id} VV     @VV \mu V \\
\mathcal{H}\otimes \mathcal{H} @ > \mu >>  \mathcal{H}
\end{CD} \qquad \qquad
\begin{CD}
\mathcal{H}\otimes \Bbbk @= \mathcal{H}  @= \Bbbk \otimes \mathcal{H}\\
@V\mathrm{id} \otimes \eta VV @V\mathrm{id} VV @VV \eta \otimes \mathrm{id} V \\
\mathcal{H}\otimes \mathcal{H} @> \mu >>  \mathcal{H} @< \mu << \mathcal{H}\otimes \mathcal{H}
\end{CD}
\end{equation}
\item a co-associative coalgebra $(\mathcal{H},\Delta,\epsilon)$ consisting of $\Bbbk$-linear comultiplication $\Delta \colon \mathcal{H} \to \mathcal{H}\otimes \mathcal{H}$ and $\Bbbk$-linear counit map $\epsilon \colon \mathcal{H} \to \Bbbk$ satisfying properties described by the diagrams:
\begin{equation}
\begin{CD}
\mathcal{H} @ > \Delta >> \mathcal{H}\otimes \mathcal{H} \\
@V \Delta  VV     @VV \Delta \otimes \mathrm{id} V \\
 \mathcal{H} \otimes \mathcal{H} @ > \mathrm{id} \otimes \Delta >> \mathcal{H}\otimes \mathcal{H}\otimes \mathcal{H} 
\end{CD} \qquad \qquad
\begin{CD}
\mathcal{H}\otimes \mathcal{H} @ < \Delta <<  \mathcal{H} @ > \Delta >> \mathcal{H}\otimes \mathcal{H} \\
@V\mathrm{id} \otimes \epsilon VV @V\mathrm{id} VV @VV \epsilon \otimes \mathrm{id} V \\
\mathcal{H}\otimes \Bbbk @= \mathcal{H}  @= \Bbbk \otimes \mathcal{H}
\end{CD}
\end{equation}
\item such that $\Delta\colon \mathcal{H} \to \mathcal{H} \otimes \mathcal{H}$ and $\epsilon \colon \mathcal{H} \to \Bbbk$ are unital algebra morphisms.
\end{enumerate}
Bialgebra $\mathcal{H}$ is \emph{graded} if it is graded as $\Bbbk$-module $\mathcal{H} = \bigoplus_{n\geq 0} \mathcal{H}^{(n)}$ with the structure maps respecting the gradation
\begin{equation}
\mathcal{H}^{(n)} \otimes \mathcal{H}^{(m)} \xrightarrow{\mu} \mathcal{H}^{(n+m)}, \qquad \mathcal{H}^{(n)} \xrightarrow{\Delta} 
\bigoplus_{n^\prime + n^{\prime\prime} = n} \mathcal{H}^{(n^\prime)} \otimes \mathcal{H}^{(n^{\prime\prime})}.
\end{equation}
A graded bialgebra is \emph{connected} if $\mathcal{H}^{(0)} \cong \Bbbk$.

The space of $\Bbbk$-linear operators $\mathrm{End} (\mathcal{H})$ can be equipped with the \emph{convolution product} 
$\star \colon \mathrm{End} (\mathcal{H}) \otimes \mathrm{End} (\mathcal{H}) \to \mathrm{End} (\mathcal{H})$ defined for $f,g \in \mathrm{End} (\mathcal{A})$
as follows
\begin{equation}
 f \star g  = \mu \circ (f \otimes g) \circ \Delta  .
\end{equation}
Such a product is associative with neutral element $\eta \circ \epsilon$. 
A bialgebra $\mathcal{H}$ is called a \emph{Hopf algebra} if there is an element $S\in \mathrm{End}_\Bbbk (\mathcal{H})$, called \emph{antipode}, which is two-sided inverse under $\star$ for the identity map $\mathrm{id}_{\mathcal{H}}$, which means
\begin{equation}
\mathrm{id}_{\mathcal{H}} \star S = S \star \mathrm{id}_{\mathcal{H}} = \eta \circ \epsilon.
\end{equation}
When it exists, the antipode $S$ is unique and is algebra anti-endomorphism: $S(1) = 1$, and $S(ab) = S(b) S(a)$ for all $a,b\in \mathcal{H}$. It is known~\cite{MilnorMoore,Takeuchi} that any graded connected bialgebra is a Hopf algebra. In that case the antipode of any homogeneous element $x \in \mathcal{H}^{(n)}$ of degree $n>0$ can be calculated recursively by 
\begin{equation} \label{eq:antipode-S}
S(x) = - x - \sum_i S(y_i)z_i = -x - \sum_i y_i S(z_i),
\end{equation}
where
\begin{equation}
\Delta x = 1\otimes x + \sum_{i} y_i \otimes z_i + x \otimes 1,
\end{equation}
and $y_i,z_i$ have degrees less then $n$.
\begin{Ex} \label{ex:FHA}
	Let $A=\{a_1, \dots , a_m \}$ be a finite set, called \emph{alphabet}, whose elements will be called \emph{letters}. A finite sequence of letters is called a \emph{word}. The set of all words on $A$ is denoted by $A^*$ (the Kleene closure operation $*$ used here shouldn't be confused with the duality sign) and turns out to be free monoid with the \emph{concatenation product} (denoted by dot $"."$ but usually omitted). The empty sequence plays the role of the neutral element of multiplication and will be denoted by $1$. Consider  free 
	algebra $\Bbbk \langle A \rangle = (\Bbbk A^*, \: . \,)$, whose linear basis consists of words, and the multiplication is given by concatenation of words, extended by linearity. 
	
	The unique compatible comultiplication and counit in $\Bbbk \langle A \rangle$ is given on letters $a_i\in A$ by 
	\begin{equation}
	\Delta(a_i) = 1 \otimes a_i + a_i \otimes 1, \qquad 
	\epsilon(a_i) = 0 , \qquad \forall a_i \in A,
	\end{equation} 
	and extended by homomorphism to words and by linearity to the whole algebra.	
	Given word $w = a_{i_1}  \dots   a_{i_n}$ we have then
	\begin{equation} \label{eq:coproduct-free}
	\Delta (a_{i_1} \dots a_{i_n} ) = \sum_{J \subset (i_1, i_2, \dots , i_n)} w_J \otimes w_{\bar{J}},
	\end{equation}
	where the multiindex $J=(j_1 , j_2 , \dots , j_k )$ is a subsequence of $ (i_1, i_2, \dots , i_n)$, 
	$w_{J} = (a_{j_1} a_{j_2}  \dots  a_{j_k})$, and $w_{\bar{J}}$ is defined analogously for the complementary subsequence $\bar{J} $. The algebra is cocommutative, graded with gradation being the length of words 
	$|a_{i_1}  \dots   a_{i_n} | = n$, locally finite and connected. The antipode on words reads 
	$S(a_{i_1}  \dots   a_{i_n} ) = (-1)^n a_{i_n}  \dots   a_{i_1} $.
\end{Ex}

Two Hopf $\Bbbk$-algebras $\mathcal{A}$, $\mathcal{B}$ are \emph{dually paired} by a map $\langle\; , \; \rangle \colon \mathcal{B} \otimes \mathcal{A} \to \Bbbk$ if
\begin{align}
\langle \mu_\mathcal{B}(b_1, b_2) , a \rangle = \langle b_1 \otimes_\mathcal{B} b_2 , \Delta_{\mathcal{A}} (a) \rangle , & \qquad
 \langle 1_{\mathcal{B}}, a \rangle = \epsilon_{\mathcal{A}} (a), \\
\langle \Delta_{\mathcal{B}} (b), a_1 \otimes_{\mathcal{A}} a_2 \rangle = \langle b_, \mu_{\mathcal{A}} (a_1, a_2) \rangle ,  & 
\qquad  \epsilon_{\mathcal{B}} (b) = \langle b , 1_{\mathcal{A}} \rangle\\
\langle S_{\mathcal{B}} (b) , a \rangle = & \langle b , S_{\mathcal{A}}(a) \rangle 
\end{align}
which is then extended to tensor products pairwise. This means that the product of $\mathcal{A}$ and coproduct of $\mathcal{B}$ are adjoint to each other under $\langle\; , \; \rangle $, and vice-versa. Likewise, the units and counits are mutually adjoint, and the antipodes are adjoint. In such case any subalgebra of $\mathcal{A}$ gives rise to the corresponding quotient algebra in $\mathcal{B}$.

When the Hopf algebra $\mathcal{H}$ is finite dimensional then the natural pairing between the $\Bbbk$-module $\mathcal{H}$ and its dual $\mathcal{H}^*$ allows to introduce on the latter the dual Hopf algebra structure. 
When $\mathcal{H}$ is infinite dimensional then there is no such general construction, which is caused by the fact that the inclusion $\mathcal{H}^* \otimes \mathcal{H}^* \subset (\mathcal{H} \otimes \mathcal{H})^*$ fails to be equality. 
For connected graded Hopf algebra $\mathcal{H} = \bigoplus_{n\geq 0} \mathcal{H}^{(n)}$ which is locally finite (each homogeneous component $\mathcal{H}^{(n)}$ is finite dimensional), one can define its \emph{graded dual} as $\mathcal{H}^{gr} = \bigoplus_{n\geq 0} \mathcal{H}^{(n)*}$ which has the property that $\mathcal{H}^{gr} \otimes \mathcal{H}^{gr} =(\mathcal{H}\otimes \mathcal{H})^{gr}$ and $(\mathcal{H}^{gr})^{gr} \cong \mathcal{H}$. Then $\mathcal{H}^{gr} \subset \mathcal{H}^*$ is a Hopf algebra where the evaluation map $\mathcal{H}^{gr} \otimes \mathcal{H} \to \Bbbk$ provides a duality pairing of $\mathcal{H}$ with $\mathcal{H}^{gr}$.
\begin{Rem}
	In Section~\ref{sec:restricted} we will consider also another construction of a dual Hopf algebra, called the \emph{restricted (or Sweedler's) dual~}\cite{Abe,Sweedler}.
\end{Rem}

\begin{Ex} \label{ex:FHA-D}
	The graded dual do the Hopf free algebra $\Bbbk \langle A \rangle$ is described as follows. By standard abuse of notation one identifies a fixed linear basis of a finite dimensional space with its dual. The dual (deconcatenation) coproduct $\delta$ is given on words by
	\begin{equation}
	\delta  (a_{i_1} \dots a_{i_n} ) = \sum_{k=0}^n  a_{i_1} \dots a_{i_k} \otimes 
	a_{i_{k+1} } \dots a_{i_n}.
	\end{equation}
	The corresponding product $\shuffle$ (called \emph{shuffle product}) dual to the coproduct $\Delta$ is given by
	\begin{equation}
	a_{i_1} \dots a_{i_{k}}  \shuffle \, 
	a_{j_1} \dots a_{j_{l}} = 
	\sum_{I = (k_1, k_2 , \dots , k_{n})} 
	a_{k_1} a_{k_2} \dots a_{k_{n}} \; ,
	\end{equation}
	where summation is over all sequences 
	$I= (k_1, k_2, \dots , k_{n} )$ such that 
	$(i_1 , i_2 , \dots , i_{k} ) \subset I$ is its subsequence, and $(j_1, j_2, \dots , j_{l}) $ is the complementary subsequence. 
	The unit, counit and antipode in the graded dual are the same as in the previous example.
\end{Ex}
\begin{Rem}
	It is known \cite{Lothaire} that the shuffle product of words can be defined recursively for all words $u,v$ and all letters $a,b$ by
	\begin{equation*}
	u\shuffle 1 = 1 \shuffle u = u, \qquad ua \shuffle vb = (ua\shuffle v)b + (u\shuffle vb)a.
	\end{equation*}
\end{Rem}
\begin{Rem}
	To distinguish between the free Hopf algebra and its graded dual, we denote them by $(\Bbbk A^*, \: . \: , \Delta)$ and $(\Bbbk A^*, \shuffle , \delta)$ respectively, skipping the unit and counit symbols.
\end{Rem}

\subsection{Hopf algebra structures on rooted ordered coloured trees} 
\label{sec:trees-Dyck}
Below we present (slightly reformulated --- see the first Remark after Proposition \ref{prop:Foissy-Hopf}) results by Foissy~\cite{Foissy-I} relevant to our paper.  

A rooted ordered tree (called also rooted plane tree) is a finite rooted tree 
$t$ such that for each vertex $v$ of $t$, the children of $v$ are totally ordered (from left to right on our pictures). Together with the depth partial order (defined by the distance from the root) this induces linear order on the vertex set $V(t)$ of the tree obtained from left-to-right depth-first search; see Figure~\ref{fig:tree-order-Dyck}. By the trivial rooted tree we understand the tree consisting of the root only. A planted rooted tree is a non-trivial rooted tree such that its root has only one child. 
\begin{figure}[h!] 
\begin{picture}(30,50)(0,0) 
\put(0,0){\line(1,2){10}}
\put(0,0){\line(-1,2){10}}
\put(-10,20){\line(-1,2){10}}
\put(10,20){\line(1,2){10}}
\put(10,20){\line(-1,2){10}}
\put(0,0){\circle*{5}}
\put(-10,-4){$1$}
\put(10,20){\circle*{5}}
\put(14,18){$4$}
\put(-20,40){\circle*{5}}
\put(-30,43){$3$}
\put(20,40){\circle*{5}}
\put(23,43){$6$}
\put(-10,20){\circle*{5}}
\put(-20,18){$2$}
\put(0,40){\circle*{5}}
\put(-7,44){$5$}
\end{picture}  
\caption{A rooted ordered tree with induced natural linear order on the vertex set}
\label{fig:tree-order-Dyck}
\end{figure}

A rooted ordered coloured (ROC) tree is a rooted tree $t$ together with a function from a set $E(t)$ of its edges to the set $A$ of colours, we assume $|A|<\infty$. By $\Bbbk T_A$ denote the linear space of finite formal combinations of A-coloured rooted ordered trees with coefficients in the field $\Bbbk$. The space $\Bbbk T_A$ is graded with the weight $|t|$ of a ROC-tree $t$ being the number of its edges
\begin{equation}
\Bbbk T_A = \oplus_{k\geq 0} \Bbbk T_A^{(k)}.
\end{equation}
By the well known connection \cite{Stanley} between rooted ordered trees and Catalan numbers $C_k$, dimension of each graded component $\Bbbk T_A^{(k)}$ equals 
\begin{equation}
\dim \Bbbk T_A^{(k)} = |A|^k C_k = \frac{|A|^k}{k+1} \left( \begin{array}{c} 2k \\ k \end{array}  \right).
\end{equation}

Define the product "$\cdot$" on $\Bbbk T_A$ as the concatenation of trees by identification of their roots; see Figure~\ref{fig:multiplication} for an example.
\begin{figure}[h!] 
\begin{picture}(20,40)(0,0) 
\color{red}
\put(8,-2){\line(0,1){20}}
\put(11,7){$b$}
\put(8,18){\line(1,2){10}}
\put(17,25){$b$}
\color{blue}
\put(8,18){\line(-1,2){10}}
\put(-5,25){$a$}
\color{black}
\put(8,-2){\circle*{5}}
\put(8,18){\circle*{5}}
\put(18,38){\circle*{5}}
\put(-2,38){\circle*{5}}
\end{picture} $.$
\begin{picture}(30,40)(-5,0) 
\color{red}
\put(8,-2){\line(-1,2){10}}
\put(-6,7){$b$}
\color{green}
\put(8,-2){\line(1,2){10}}
\put(17,7){$c$}
\color{black}
\put(8,-2){\circle*{5}}
\put(-2,18){\circle*{5}}
\put(18,18){\circle*{5}}
\end{picture}    $= $
\begin{picture}(70,40)(-30,0) 
\color{red}
\put(8,-2){\line(-1,1){20}}
\put(-12,18){\line(1,2){10}}
\put(-13,6){$b$}
\put(2,7){$b$}
\color{blue}
\put(-12,18){\line(-1,2){10}}
\color{red}
\put(8,-2){\line(0,1){20}}
\put(-2,25){$b$}
\color{blue}
\put(-27,25){$a$}
\color{green}
\put(8,-2){\line(1,1){20}}
\put(23,7){$c$}
\color{black}
\put(8,-2){\circle*{5}}
\put(-12,18){\circle*{5}}
\put(-22,38){\circle*{5}}
\put(8,18){\circle*{5}}
\put(-2,38){\circle*{5}}
\put(28,18){\circle*{5}}
\end{picture} 
\caption{Multiplication of two coloured ordered rooted trees}
\label{fig:multiplication}
\end{figure}
The product respects the gradation,
%\begin{equation}
%\Bbbk T_n^{(k_1)} . \; \Bbbk T_n^{(k_2)} \subset \Bbbk T_n^{(k_1+k_2)} ,
%\end{equation}
is associative with the trivial tree \ttree being the neutral element (i.e. the unit map $\eta \colon \Bbbk \to \Bbbk T_A$ is defined by $1 \longmapsto \ttree$).

\begin{figure}[h!] 
\begin{picture}(70,40)(-30,0) 
\color{red}
\put(8,-2){\line(-1,1){20}}
\put(-12,18){\line(1,2){10}}
\put(-13,6){$b$}
\put(2,7){$b$}
\color{blue}
\put(-11.5,18){\line(-1,2){10}}
\put(-12.5,18){\line(-1,2){10}}
\color{red}
\put(7.5,-2){\line(0,1){20}}
\put(8.5,-2){\line(0,1){20}}
\put(-2,25){$b$}
\color{blue}
\put(-27,25){$a$}
\color{green}
\put(8,-2){\line(1,1){20}}
\put(23,7){$c$}
\color{black}
\put(8,-2){\circle*{5}}
\put(-12,18){\circle*{5}}
\put(-22,38){\circle*{5}}
\put(8,18){\circle*{5}}
\put(-2,38){\circle*{5}}
\put(28,18){\circle*{5}}
\end{picture}    $\xrightarrow{\text{pruning}} \qquad $
\begin{picture}(25,40)(0,0) 
\color{red}
\put(8,-2){\line(1,2){10}}
\put(17,5){$b$}
\color{blue}
\put(8,-2){\line(-1,2){10}}
\put(-5,5){$a$}
\color{black}
\put(8,-2){\circle*{5}}
\put(18,18){\circle*{5}}
\put(-2,18){\circle*{5}}
\end{picture} $\otimes$
\begin{picture}(30,40)(-10,0) 
\color{red}
\put(8,-2){\line(-1,2){10}}
\put(-6,5){$b$}
\put(-2,18){\line(0,1){20}}
\put(-10,25){$b$}
\color{green}
\put(8,-2){\line(1,2){10}}
\put(17,5){$c$}
\color{black}
\put(8,-2){\circle*{5}}
\put(-2,18){\circle*{5}}
\put(-2,38){\circle*{5}}
\put(18,18){\circle*{5}}
\end{picture} 
\caption{Pruning of a ROC-tree; pruned branches are thickened}
\label{fig:pruning}
\end{figure}
In order to define compatible coproduct on $\Bbbk T_A$ one has first to describe the operation of pruning of a tree. A rooted subtree $t_s$ of a ROC-tree $t$ is called \emph{admissible} if it shares the root of $t$. Such an admissible subtree is again ROC-tree with the root, order and colours inherited from $t$. The set of admissible subtrees of $t$ (including the trivial tree and $t$ itself) will be denoted by $A(t)$. Given such admissible subtree $t_s\subset t$ it defines a sequence $(t_1,\dots ,t_m)$ of planted trees being branches of $t$ pruned to get $t_s$, with the order in the sequence inherited from the order on $t$. By concatenation of the pruned branches we obtain the \emph{complementary tree} $t_c = t_1 \cdot \ldots \cdot  t_m$ to the admissible subtree $t_s$ of $t$. Such a pruning operation gives an element $t_c\otimes t_s$; see Figure~\ref{fig:pruning} for an example.

The coproduct of a tree is defined as sum of pairs $t_c \otimes t_s$ for all admissible subtrees of $t$; see Figure~\ref{fig:coproduct}
\begin{equation} \label{eq:Delta-cop}
\Delta(t) = \sum_{t_s\in A(t)} t_c \otimes t_s ,
\end{equation}
\begin{figure}[h!] 
\begin{equation*}
\Delta \left( 
\begin{picture}(30,40)(2,10)
\color{red}
\put(18,0){\line(-1,2){10}}
\put(18,0){\line(1,2){10}}
\put(3,8){$b$}
\put(28,8){$b$}
\color{blue}
\put(8,20){\line(0,1){20}}
\put(12,29){$a$}
\color{black}
\put(8,20){\circle*{5}}
\put(8,40){\circle*{5}}
\put(18,0){\circle*{5}}
\put(28,20){\circle*{5}}
\end{picture} 
\right)  
=  
\begin{picture}(32,40)(2,10)
\color{red}
\put(18,0){\line(-1,2){10}}
\put(18,0){\line(1,2){10}}
\put(3,8){$b$}
\put(28,8){$b$}
\color{blue}
\put(8,20){\line(0,1){20}}
\put(12,29){$a$}
\color{black}
\put(8,20){\circle*{5}}
\put(8,40){\circle*{5}}
\put(18,0){\circle*{5}}
\put(28,20){\circle*{5}}
\end{picture} 
\otimes
\begin{picture}(12,3)(3,0)
\put(8,2){\circle*{5}}
\end{picture} 
+
\begin{picture}(12,40)(5,10)
\color{red}
\put(8,0){\line(0,1){20}}
\put(12,8){$b$}
\color{blue}
\put(8,20){\line(0,1){20}}
\put(12,28){$a$}
\color{black}
\put(8,20){\circle*{5}}
\put(8,40){\circle*{5}}
\put(8,0){\circle*{5}}
\end{picture} 
\otimes
\begin{picture}(12,40)(5,10)
\color{red}
\put(8,0){\line(0,1){20}}
\put(12,8){$b$}
\color{black}
\put(8,20){\circle*{5}}
\put(8,0){\circle*{5}}
\end{picture} 
+
\begin{picture}(33,40)(2,10)
\color{red}
\put(18,0){\line(1,2){10}}
\put(28,8){$b$}
\color{blue}
\put(18,0){\line(-1,2){10}}
\put(3,8){$a$}
\color{black}
\put(8,20){\circle*{5}}
\put(18,0){\circle*{5}}
\put(28,20){\circle*{5}}
\end{picture} 
\otimes
\begin{picture}(12,40)(5,10)
\color{red}
\put(8,0){\line(0,1){20}}
\put(12,8){$b$}
\color{black}
\put(8,20){\circle*{5}}
\put(8,0){\circle*{5}}
\end{picture} 
+
\begin{picture}(12,40)(5,10)
\color{blue}
\put(8,0){\line(0,1){20}}
\put(12,8){$a$}
\color{black}
\put(8,20){\circle*{5}}
\put(8,0){\circle*{5}}
\end{picture} 
\otimes
\begin{picture}(33,40)(2,10)
\color{red}
\put(18,0){\line(-1,2){10}}
\put(18,0){\line(1,2){10}}
\put(3,8){$b$}
\put(28,8){$b$}
\color{black}
\put(8,20){\circle*{5}}
\put(18,0){\circle*{5}}
\put(28,20){\circle*{5}}
\end{picture} 
+
\begin{picture}(12,40)(5,10)
\color{red}
\put(8,0){\line(0,1){20}}
\put(12,8){$b$}
\color{black}
\put(8,20){\circle*{5}}
\put(8,0){\circle*{5}}
\end{picture} 
\otimes
\begin{picture}(12,40)(5,10)
\color{red}
\put(8,0){\line(0,1){20}}
\put(12,8){$b$}
\color{blue}
\put(8,20){\line(0,1){20}}
\put(12,28){$a$}
\color{black}
\put(8,20){\circle*{5}}
\put(8,40){\circle*{5}}
\put(8,0){\circle*{5}}
\end{picture} 
+
\begin{picture}(10,3)(2,0)
\put(8,2){\circle*{5}}
\end{picture} 
\otimes
\begin{picture}(30,40)(2,10)
\color{red}
\put(18,0){\line(-1,2){10}}
\put(18,0){\line(1,2){10}}
\put(3,8){$b$}
\put(28,8){$b$}
\color{blue}
\put(8,20){\line(0,1){20}}
\put(12,29){$a$}
\color{black}
\put(8,20){\circle*{5}}
\put(8,40){\circle*{5}}
\put(18,0){\circle*{5}}
\put(28,20){\circle*{5}}
\end{picture} 
\end{equation*}
\caption{The pruning coproduct of a ROC-tree}
\label{fig:coproduct}
\end{figure}
and then extended to $\Bbbk T_A$ by linearity.

Such coproduct is coassociative, respects the gradation,
%\begin{equation}
%\Delta( \Bbbk T_n^{(k)}) \subset \bigoplus_{k_1 + k_2 = k} \Bbbk T_n^{(k_1)} \otimes %\Bbbk T_n^{(k_2)} ,
%\end{equation}
and is compatible with the counit defined on trees as
\begin{equation} \label{eq:counit}
\epsilon(t) = \begin{cases} 1 & \text{if} \quad t = \ttree , \\
0 & \text{otherwise.} \end{cases}
\end{equation}
In this context it is convenient to define the operation $B^+_i$ of planting of a tree on a new root by attaching it to the old one by additional edge coloured by $i$. In particular, planting allows to define the coproduct recursively starting from $\Delta(\ttree) = \ttree \otimes \ttree$, and using then the formula
\begin{equation} \label{eq:B+}
\Delta(B_i^+(t)) = B_i^+ (t) \otimes \ttree + (\mathrm{id} \otimes B_i^+) \circ \Delta(t),
\end{equation}
together with compatibility of the coproduct $\Delta$ with the concatenation product. Equation \eqref{eq:B+} has the following simple meaning: apart from the trivial subtree, all admissible subtrees of a planted tree contain the lowest (i.e. incident to the root) edge.
\begin{Prop} \label{prop:Foissy-Hopf}
The concatenation multiplication and pruning coproduct with the corresponding unit and counit maps equip $\Bbbk T_A$ with the structure of graded locally finite and connected bialgebra (thus Hopf algebra).
\end{Prop}
\begin{Rem}
	The above result was given by Foissy~\cite{Foissy-I} in the equivalent setting of  the rooted ordered vertex-coloured (or decorated) forests.  Any ROC tree is uniquely mapped, by deletion of the root, to an ordered forest colouring first its vertices using colours of adjacent edges below them; see the bijection map visualized on 
	Figure~\ref{fig:bijections}. This notation resulted as decorated and non-commutative version of the Connes--Kreimer Hopf algebra~\cite{Connes-Kreimer} used to explain the renormalization procedure in the quantum field theory.
\end{Rem}
\begin{figure}[h!] 
\begin{picture}(70,60)(-30,0) 
\color{red}
\put(8,-2){\line(-1,1){20}}
\put(-12,18){\line(1,2){10}}
\put(-13,6){$b$}
\put(2,7){$b$}
\color{blue}
\put(-12,18){\line(-1,2){10}}
\color{green}
\put(-2,38){\line(0,1){20}}
\put(23,6){$c$}
\color{red}
\put(8,-2){\line(0,1){20}}
\put(-2,25){$b$}
\color{blue}
\put(28,18){\line(0,1){20}}
\put(32,25){$a$}
\put(-27,25){$a$}
\color{green}
\put(8,-2){\line(1,1){20}}
\put(2,47){$c$}
\color{black}
\put(8,-2){\circle*{5}}
\put(-12,18){\circle*{5}}
\put(28,38){\circle*{5}}
\put(-22,38){\circle*{5}}
\put(-2,58){\circle*{5}}
\put(8,18){\circle*{5}}
\put(-2,38){\circle*{5}}
\put(28,18){\circle*{5}}
\end{picture} $\longleftrightarrow \; $
\begin{picture}(70,50)(-30,20) 
\put(-12,18){\line(1,2){10}}
\put(-12,18){\line(-1,2){10}}
\put(-2,38){\line(0,1){20}}
\put(28,18){\line(0,1){20}}
\color{red}
\put(-12,18){\circle*{5}}
\put(8,18){\circle*{5}}
\put(-2,38){\circle*{5}}
\put(-23,18){$b$}
\put(6,24){$b$}
\put(3,37){$b$}
\color{blue}
\put(28,38){\circle*{5}}
\put(-22,38){\circle*{5}}
\put(28,43){$a$}
\put(-23,43){$a$}
\color{green}
\put(28,18){\circle*{5}}
\put(-2,58){\circle*{5}}
\put(33,18){$c$}
\put(-2,65){$c$}
\end{picture}
\caption{Transition from the setting of ROC-trees to the setting of ROD-forests}
\label{fig:bijections}
\end{figure}
\begin{Rem}
	In \cite{Foissy-I} one can find, among others, also the corresponding description of the antipode, which can be transferred from the ROD-forests to ROC-trees. 
\end{Rem}
\begin{Cor} \label{cor:kA-t}
	The subalgebra of $\Bbbk T_A$ generated by one-edge planted trees 
	$a_i = B_i^+ (\ttree)$ is a Hopf subalgebra isomorphic to the free Hopf algebra described in Example~\ref{ex:FHA}.
\end{Cor}

We conclude this Section by presenting the graded dual of the Hopf algebra of ROC-trees, which is again reformulation of the corresponding results of~\cite{Foissy-I}. Because the natural basis of the finite dimensional subspace $\Bbbk T_A^{(k)}$ is provided by ROC-trees of weight $k$ it seems natural to represent the dual basis of $(\Bbbk T_A^{(k)})^*$ by such trees again, i.e. the functional
$\phi_t \in (\Bbbk T_A^{(k)})^*$ defined on trees by
\begin{equation}
\langle \phi_t , t^\prime \rangle = \begin{cases} 1 & \text{if} \quad t = t^\prime , \\
0 & \text{otherwise} \end{cases}
\end{equation}
by standard abuse of notation is identified with $t$. The dual (deconcatenation) coproduct $\delta=(.)^*$ to the concatenation product acts on trees as
\begin{equation} \label{eq:delta-cop}
\delta(t) = \sum_{t^\prime , t^{\prime\prime} \in T_n} 
\langle t , t^\prime . \; t^{\prime\prime} \rangle \; t^\prime \otimes t^{\prime\prime} =
\sum_{t^\prime \cdot \; t^{\prime\prime} = t} 
t^\prime \otimes t^{\prime\prime},
\end{equation}
i.e. when $t=t_1\dots t_m$ is planted trees decomposition, then
\begin{equation}
\delta( t_1\dots t_m ) = \sum_{i=0}^m (t_1\dots t_i) \otimes 
(t_{i+1}\dots t_m) ,
\end{equation}
see Figure~\ref{fig:coproduct-delta}.
\begin{figure}[h!]
\begin{align*}
\delta\left( 
\begin{picture}(30,40)(2,10)
\color{red}
\put(18,0){\line(-1,2){10}}
\put(18,0){\line(1,2){10}}
\put(3,8){$b$}
\put(28,8){$b$}
\color{blue}
\put(8,20){\line(0,1){20}}
\put(12,29){$a$}
\color{black}
\put(8,20){\circle*{5}}
\put(8,40){\circle*{5}}
\put(18,0){\circle*{5}}
\put(28,20){\circle*{5}}
\end{picture} 
\right) & 
= 
\begin{picture}(10,3)(2,0)
\put(8,2){\circle*{5}}
\end{picture} 
\otimes
\begin{picture}(30,40)(2,10)
\color{red}
\put(18,0){\line(-1,2){10}}
\put(18,0){\line(1,2){10}}
\put(3,8){$b$}
\put(28,8){$b$}
\color{blue}
\put(8,20){\line(0,1){20}}
\put(12,29){$a$}
\color{black}
\put(8,20){\circle*{5}}
\put(8,40){\circle*{5}}
\put(18,0){\circle*{5}}
\put(28,20){\circle*{5}}
\end{picture} + 
\begin{picture}(12,40)(5,10)
\put(8,20){\circle*{5}}
\put(8,40){\circle*{5}}
\put(8,0){\circle*{5}}
\color{red}
\put(8,0){\line(0,1){20}}
\put(12,8){$b$}
\color{blue}
\put(8,20){\line(0,1){20}}
\put(12,28){$a$}
\end{picture} 
\otimes
\begin{picture}(12,40)(5,10)
\put(8,20){\circle*{5}}
\put(8,0){\circle*{5}}
\color{red}
\put(8,0){\line(0,1){20}}
\put(12,8){$b$}
\end{picture} 
 +
\begin{picture}(32,40)(2,10)
\color{red}
\put(18,0){\line(-1,2){10}}
\put(18,0){\line(1,2){10}}
\put(3,8){$b$}
\put(28,8){$b$}
\color{blue}
\put(8,20){\line(0,1){20}}
\put(12,29){$a$}
\color{black}
\put(8,20){\circle*{5}}
\put(8,40){\circle*{5}}
\put(18,0){\circle*{5}}
\put(28,20){\circle*{5}}
\end{picture} 
\otimes
\begin{picture}(12,3)(3,0)
\put(8,2){\circle*{5}}
\end{picture} \end{align*}
\caption{The deconcatenation coproduct of a ROC-tree}
\label{fig:coproduct-delta}
\end{figure}
\begin{Cor} \label{cor:delta-.}
	Equation \eqref{eq:delta-cop} implies the following matching condition between the deconcatenation coproduct and the concatenation product
	\begin{equation}
	\delta(s.t) = \delta(s).(\ttree \otimes t) + (s\otimes \ttree). \delta(t) - s\otimes t, \qquad s,t\in T_A,
	\end{equation}
where, by the standard abuse of notation, we extended the product sign from $\Bbbk T_A$ to $\Bbbk T_A \otimes \Bbbk T_A$.
\end{Cor}

The (asymmetric shuffle or grafting) product $\Tshuffle = \Delta^*$, dual to the pruning coproduct satisfies
\begin{equation} \label{eq:ashuffle}
\Delta(t) = \sum_{t^\prime , t^{\prime\prime} \in T_A} 
\langle t , t^\prime \Tshuffle t^{\prime\prime} \rangle \; t^\prime \otimes t^{\prime\prime} ,
\end{equation}
and is defined with the help of the grafting procedure that follows from comparison of equations \eqref{eq:Delta-cop} and \eqref{eq:ashuffle}. Given ROC-tree $t^\prime=t_1 . \dots , t_m$ decomposed into the planted factors, its (non-unique) grafting on ROC-tree $t^{\prime\prime}$ is defined as attaching roots of the factors $t_i$ to vertices of $t^{\prime\prime}$ in a way, which preserves the original ordering of the factors. In other words, a grafting of $t^\prime$ on $t^{\prime\prime}$ gives a tree $\tilde{t}$ such that there exists a pruning with $t^{\prime\prime} = \tilde{t}_s$ with the corresponding $t^\prime = \tilde{t}_c$; see Figure~\ref{fig:ashuffle-prod} for an example. 
\begin{figure}[h!]
\begin{picture}(30,50)(-30,0) 
\color{red}
\put(-12,0){\line(1,2){10}}
\put(-4,8){$b$}
\color{blue}
\put(-12,0){\line(-1,2){10}}
\put(-26,8){$a$}
\color{black}
\put(-12,0){\circle*{5}}
\put(-2,20){\circle*{5}}
\put(-22,20){\circle*{5}}
\end{picture}
$\Tshuffle$
\begin{picture}(20,40)
\color{red}
\put(8,0){\line(0,1){20}}
\put(12,8){$b$}
\color{black}
\put(8,0){\circle*{5}}
\put(8,20){\circle*{5}}
\end{picture}
$=$
\begin{picture}(50,40)
\color{red}
\put(27.5,0){\line(0,1){20}}
\put(28.5,0){\line(0,1){20}}
\put(31,9){$b$}
\put(28,0){\line(1,1){20}}
\put(44,7){$b$}
\color{blue}
\put(27.5,0){\line(-1,1){20}}
\put(28.5,0){\line(-1,1){20}}
\put(8,7){$a$}
\color{black}
\put(28,0){\circle*{5}}
\put(28,20){\circle*{5}}
\put(8,20){\circle*{5}}
\put(48,20){\circle*{5}}
\end{picture}
$+$
\begin{picture}(30,40)
\color{red}
\put(18,0){\line(1,2){10}}
\put(27.5,20){\line(0,1){20}}
\put(28.5,20){\line(0,1){20}}
\put(28,8){$b$}
\put(32,29){$b$}
\color{blue}
\put(17.5,0){\line(-1,2){10}}
\put(18.5,0){\line(-1,2){10}}
\put(3,8){$a$}
\color{black}
\put(8,20){\circle*{5}}
\put(28,40){\circle*{5}}
\put(18,0){\circle*{5}}
\put(28,20){\circle*{5}}
\end{picture}
$+$
\begin{picture}(50,40)
\color{red}
\put(28,0){\line(0,1){20}}
\put(31,9){$b$}
\put(27.5,0){\line(1,1){20}}
\put(28.5,0){\line(1,1){20}}
\put(44,7){$b$}
\color{blue}
\put(27.5,0){\line(-1,1){20}}
\put(28.5,0){\line(-1,1){20}}
\put(8,7){$a$}
\color{black}
\put(28,0){\circle*{5}}
\put(28,20){\circle*{5}}
\put(8,20){\circle*{5}}
\put(48,20){\circle*{5}}
\end{picture}
$+$
\begin{picture}(20,40)
\color{red}
\put(8,0){\line(0,1){20}}
\put(7.5,20){\line(1,2){10}}
\put(8.5,20){\line(1,2){10}}
\put(12,8){$b$}
\put(17,28){$b$}
\color{blue}
\put(7.5,20){\line(-1,2){10}}
\put(8.5,20){\line(-1,2){10}}
\put(-6,28){$a$}
\color{black}
\put(8,0){\circle*{5}}
\put(8,20){\circle*{5}}
\put(18,40){\circle*{5}}
\put(-2,40){\circle*{5}}
\end{picture}
$+$
\begin{picture}(30,40)
\color{red}
\put(18,0){\line(-1,2){10}}
\put(17.5,0){\line(1,2){10}}
\put(18.5,0){\line(1,2){10}}
\put(3,8){$b$}
\put(28,8){$b$}
\color{blue}
\put(7.5,20){\line(0,1){20}}
\put(8.5,20){\line(0,1){20}}
\put(12,29){$a$}
\color{black}
\put(8,20){\circle*{5}}
\put(8,40){\circle*{5}}
\put(18,0){\circle*{5}}
\put(28,20){\circle*{5}}
\end{picture}
$+$
\begin{picture}(40,40)
\color{red}
\put(28,0){\line(-1,1){20}}
\put(27.5,0){\line(1,1){20}}
\put(28.5,0){\line(1,1){20}}
\put(8,7){$b$}
\put(44,7){$b$}
\color{blue}
\put(27.5,0){\line(0,1){20}}
\put(28.5,0){\line(0,1){20}}
\put(31,9){$a$}
\color{black}
\put(28,0){\circle*{5}}
\put(28,20){\circle*{5}}
\put(8,20){\circle*{5}}
\put(48,20){\circle*{5}}
\end{picture}
\caption{The asymmetric shuffle (or grafting) product of two ROC-trees; grafted branches are thickened}
\label{fig:ashuffle-prod}
\end{figure}
With the grafting product $\Tshuffle$, deconcatenation coproduct $\delta$, the unit $\eta = \epsilon^*$ and counit $\epsilon = \eta^*$ maps, the space spanned by ROC-trees is equipped with another bialgebra (thus Hopf algebra) structure -- the graded dual to the previous one.
\begin{Rem}
It is remarkable fact, discovered by Foissy~\cite{Foissy-I}, that the duality described above is self-duality. The situation is analogous to the well known self-duality of the Hopf algebra of symmetric functions~\cite{Stanley}. 
\end{Rem}

\subsection{Hopf algebra of quasi-symmetric functions, and its graded dual} \label{sec:qs-ns}
Let $x= (x_1, x_2, x_3, \dots )$ denote infinite \emph{totally ordered} set of commuting variables, and let $\Bbbk [[x_1, x_2, x_3, \dots ]]$ be the algebra of formal power series of bounded degree. Such a formal series is called \emph{quasi-symmetric function} if the coefficient of any term $x_{i_1}^{n_1} x_{i_2}^{n_2} \dots x_{i_k}^{n_k}$ 
with $i_1 < i_2 < \dots < i_k$ \emph{strictly increasing}, agrees with that of $x_1^{n_1} x_2^{n_2} \dots x_k^{n_k}$. The linear space $\mathrm{QSym}$ of quasi-symmetric functions has as a basis the \emph{monomial quasi-symmetric functions} indexed by compositions. Recall that a composition of $m$, written $\alpha \models m$, is a finite sequences of positive integers $\alpha=(\alpha_1, \alpha_2 , \dots , \alpha_k)$ such that $|\alpha| = \alpha_1 + \alpha_2 + \dots + \alpha_k = m$. In this case we say that $\alpha$ has $k$ parts, or it is of length $\ell(\alpha)=k$. The elements of the basis are of the form
\begin{equation}
M_\alpha = \sum_{i_1 < i_2 < \dots < i_k} 
x_{i_1}^{\alpha_1} x_{i_2}^{\alpha_2} \dots x_{i_k}^{\alpha_k},
\end{equation}
where the sum is over all $k$-tuples $(i_1, i_2, \dots , i_k)$ of strictly increasing; by definition $M_\emptyset = 1$. The algebra $\mathrm{QSym}$ is graded with each graded component $\mathrm{QSym}^{(m)}$ spanned by those $M_\alpha$ for which 
$|\alpha|=m$. By the well known bijection~\cite{Introduction-QSym} any such composition can be identified with a subset 
\begin{equation}
\mathrm{set}(\alpha) = 
\{\alpha_1,  \alpha_1 + \alpha_2, \dots ,  \alpha_1 + \dots + \alpha_{k-1} \}
\end{equation} 
of $\{1,2,\dots ,m-1\}$, therefore $\dim \mathrm{QSym}^{(m)} = 2^{m-1}$.

We can put a partial order on the set of all compositions of $m$ by refinement. The covering relations are of the form 
\begin{equation}
(\alpha_1, \dots , \alpha_i, \alpha_{i+1}, \dots \alpha_k) \prec (\alpha_1, \dots , \alpha_i + \alpha_{i+1}, \dots \alpha_k).
\end{equation}
This allows to define another important basis formed by the \emph{fundamental} quasi-symmetric functions, also indexed by compositions
\begin{equation}
F_\alpha = \sum_{\beta \preccurlyeq \alpha} M_\beta.
\end{equation}
By inclusion-exclusion we can express the $M_\alpha$ in terms of the $F_\alpha$
\begin{equation}
M_\alpha = \sum_{\beta \preccurlyeq \alpha} (-1)^{\ell(\beta) - \ell(\alpha)} F_\beta .
\end{equation}
 
The product in $\mathrm{QSym}$, inherited from the standard multiplication of power series, can be described in the basis $(M_\alpha)$ in terms of the quasi-shuffle (or overlapping shuffle) $\Qshuffle$ of compositions: in addition to shuffling components $\alpha_i$ and $\beta_j$ of two compositions $\alpha=(\alpha_1, \dots , \alpha_k)$ and $\beta=(\beta_1, \dots , \beta_l)$ we may replace any number of pairs of consecutive components $\alpha_i$ and $\beta_j$ in the shuffle by their sum $\alpha_i + \beta_j$
\begin{equation} \label{eq:M-a-M-b}
M_\alpha M_\beta = \sum_{\gamma} M_\gamma \;,
\end{equation}
where $\gamma$ is a summand in quasi-shuffle of $\alpha$ and $\beta$.
\begin{Ex}
For $M_{(1)} = x_1 + x_2 + \dots $ and $M_{(2)} = x_1^2 + x_2^2 + \dots $ we have
\begin{equation*}
(x_1 + x_2 + \dots ) ( x_1^2 + x_2^2 + \dots ) =
(x_1 x_2^2 + x_1 x_3^2 + \dots ) + (x_1^2 x_2 + x_1^2 x_3 + \dots ) + (x_1^3 + x_2^3 + \dots ),
\end{equation*}
therefore we obtain
\begin{equation}
M_{(1)} M_{(2)} = M_{(1,2)} + M_{(2,1)} + M_{(3)} , \qquad \text{or} \quad (1)\Qshuffle (2) = (1,2) + (2,1) + (3)\; .
\end{equation}
\end{Ex}

The coproduct $\delta$ in the algebra of quasi-symmetric functions can be defined using the doubling variables trick. Here to the totally ordered set of variables $x= ( x_1, x_2, x_3, \dots )$ we add its copy $y= ( y_1, y_2, y_3, \dots )$ placing elements of $y$ \emph{after} elements of $x$, and getting the ordered sum of the sets of variables. To obtain the coproduct $\delta(f)$ of a quasi-symmetric function $f$ we expand the function over the doubled variables, decompose resulting expression into sum of products of functions of $x$ and $y$ getting this way
\begin{equation} \label{eq:var-doubling}
f \mapsto f(x ) \mapsto f(x,y) =
\sum_j f^\prime _j (x ) f_j ^{\prime\prime} (y) \mapsto \sum_j f_j^\prime \otimes f_j ^{\prime\prime} = \delta(f).
\end{equation}
In the basis of monomial quasi-symmetric functions $(M_\alpha)$ the coproduct formula reads
\begin{equation} \label{eq:Delta-M}
\delta(M_\alpha) = \sum_{\beta \cdot \gamma = \alpha} M_\beta \otimes M_\gamma \; ,
\end{equation}
where $\beta \cdot \gamma$ is concatenation of two compositions. As a result we obtain graded, locally finite and connected bialgebra (thus Hopf algebra) which is commutative but not cocommutative.
\begin{Ex}
Applying the procedure to $M_{(2,1)} = x_1^2 x_2 + x_1^2 x_3 +  \dots $ we have
\begin{equation*} \begin{split}
x_1^2 x_2 + x_1^2 x_3 +  \dots & \mapsto  
x_1^2 x_2 + x_1^2 x_3 + \dots +  x_1^2 y_1 + x_1^2 y_2 +  \dots + y_1^2 y_2 + y_1^2 y_3 + \dots = \\
& = M_{(2,1)}(x) +   M_{(2)}(x) M_{(1)}(y)  +
M_{(2,1)}(y) 
\end{split}
\end{equation*}
getting this way 
\begin{equation}
\delta(M_{(2,1)}) = M_{(2,1)} \otimes 1 + M_{(2)} \otimes M_{(1)} + 
 1 \otimes M_{(2,1)} \; .
\end{equation}
\end{Ex}

The graded dual to $\mathrm{QSym}$ is called the Hopf algebra of non-commutative symmetric functions~\cite{Gelfand-NSym} and denoted by $\mathrm{NSym}$. Let $(H_\alpha)$ be the dual basis to $(M_\beta)$
\begin{equation}
\langle H_\alpha , M_\beta \rangle = \delta_{\alpha \beta},
\end{equation}
then by dualization of equations \eqref{eq:M-a-M-b} and \eqref{eq:Delta-M} we obtain the product and coproduct formulas in $\mathrm{NSym}$
\begin{equation} \label{eq:H-prod-copr}
H_\alpha H_\beta = H_{\alpha \cdot \beta}, \qquad
\Delta(H_\alpha ) = \sum_{(\beta, \gamma )}
 H_\beta \otimes H_\gamma \; ,
\end{equation} 
where $\alpha$ can be obtained as a summand in quasi-shuffle of $\beta$ and $\gamma$.
In particular, for a composition $\alpha = (\alpha_1, \dots , \alpha_k)$ one has
\begin{equation}
H_\alpha = H_{\alpha_1} \dots H_{\alpha_k},
\end{equation}
where we wrote $H_j = H_{(j)}$ for a composition $(j)$. The dual element to 
$M_\emptyset = 1$ is $H_{\emptyset} = H_0 = 1$. This leads to the conclusion that, as algebra, $\mathrm{NSym}$ is freely generated by non-commuting elements $H_1, H_2, \dots $
\begin{equation}
\mathrm{NSym} = \Bbbk \langle H_1, H_2, \dots \rangle .
\end{equation}
The coproduct formula for the generators follows from equations \eqref{eq:H-prod-copr} and reads
\begin{equation} \label{eq:cop-H-n}
\Delta(H_m) = \sum_{j=0}^m H_j \otimes H_{m-j}.
\end{equation}
\begin{Rem}
It is known \cite{Hoffman-nc} that the Hopf algebra $\mathrm{NSym}$ is isomorphic to a Hopf subalgebra of rooted ordered (monochromatic) trees generated by $ (B^+)^m (\ttree) \leftrightarrow H_m $, where in the monochromatic $|A|=1$ case we skip the lower index $i=1$ describing the colour of the attached edge.
\end{Rem}
\begin{Rem}
	One can recapitulate this Section in the spirit of Examples~\ref{ex:FHA} and \ref{ex:FHA-D} that we presented two, mutually dual, Hopf algebra structures in the space of compositions.
	Similarly to the shuffle product of words, the quasi-shuffle of compositions can be defined recursively~\cite{Hoffman-sh} for all compositions $\alpha,\beta$ and all natural numbers $k,l$ by
	\begin{equation*}
	\alpha\Qshuffle \emptyset = \emptyset \Qshuffle \alpha = \alpha, \quad ((k)\cdot\alpha) \Qshuffle ((l)\cdot\beta) = (k)\cdot(\alpha\Qshuffle ((l)\cdot\beta)) + (l)\cdot(((k)\cdot\alpha)\Qshuffle \beta) + (k+l)\cdot(\alpha\Qshuffle \beta).
	\end{equation*}
\end{Rem}

\section{The Hopf algebra of coloured non-commutative symmetric functions} \label{sec:n-NSyn-QSym}

\subsection{The Hopf algebra of sentences} \label{sec:sent-A}
Given finite set $A=\{ a_1, \dots , a_n \}$ called alphabet, and whose elements are called letters. Words are finite sequences of letters (written without separation), sentences are finite sequences of words (instead of spacing to separate them we use commas). The size of a word $w = a_{i_1} \dots a_{i_k}$ is the number $|w|=k$ of its letters, the size of a sentence $I = (w_1,w_2,\dots , w_m)$ is the sum $|I| = |w_1| + |w_2| + \dots +|w_m|$ of sizes of its words, while the length of the sentence is the number $\ell(I) = m$ of its words. The maximal word of the sentence is the concatenation $w(I) = w_1 w_2 \dots w_m$ of all its words. 
\begin{Ex}
The sentence $I = (aba,ca,bac)$ over alphabet $\{a,b,c\}$ has length $\ell(I) = 3$, is of size $|I|=8$. Its maximal word is $w(I)= abacabac$.	
\end{Ex}
\begin{Rem}
	For unary alphabet $A = \{a\}$ we will identify words with their size, and then identify sentences with corresponding compositions, for example $(aaa,aa,aaa) \longleftrightarrow (3,2,3)$. 
Speaking about colours in the place of letters, instead of sentences we may use the notion of coloured compositions.
\end{Rem}

Given two sentences $I$, $J$, we say that $I$ is coarsening of $J$ (or equivalently, $J$ is refinement of $I$), denoted by $I \succcurlyeq J$, if we can obtain the words of $I$ by concatenation of adjacent words of $J$. For example $(aba,ca,bac) \succcurlyeq (ab,a,ca,ba,c)$. With the refinement order, the poset of sentences having the same maximal word $w$ is isomorphic to the poset of compositions of $|w|$, and therefore isomorphic (recall the bijection $\mathrm{set}$ mentioned in Section~\ref{sec:qs-ns}) to the (dual of the) Boolean poset of $\{ 1,2 ,\dots ,|w|-1 \}$.
\begin{figure}[h!] 
\begin{center}
	\includegraphics[width=14cm]{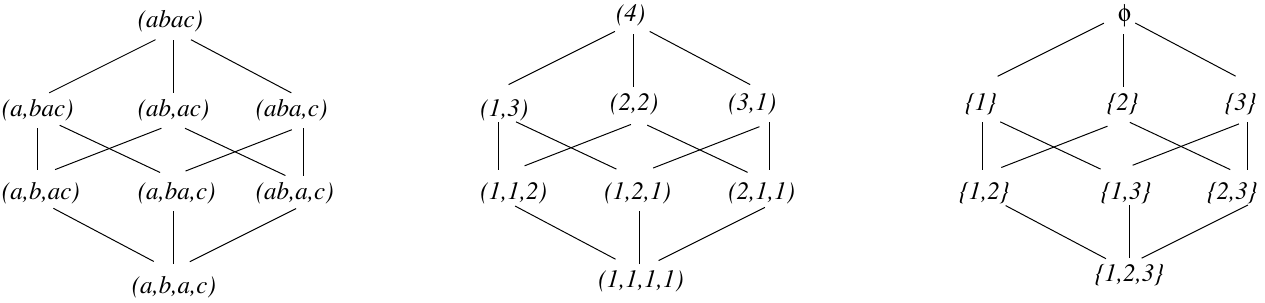}
\end{center}
\caption{Isomorphism of three posets}
\label{fig:posets}
\end{figure}
This allows to transfer structural results between the posets, in particular to find description of the Moebius function of the posets of sentences
\begin{equation} \label{eq:Moebius}
\mu(J,I) = (-1)^{\ell(J) - \ell(I)} \qquad \text{for} \quad J\preccurlyeq I.
\end{equation}
\begin{Rem}
	The above poset point of view makes the theory of coloured non-commutative symmetric and quasi-symmetric functions so similar to the original monochromatic theory. There are however various differences  related to new degree of non-commutativity caused by presence of different colours/letters. 
\end{Rem}
\begin{Cor} \label{cor:maximal-word}
	The set of sentences of size $n\in\mathbb{N}$ splits onto $|A|^n$ pairwise disjoint classes of sentences with the same maximal word.
\end{Cor}

Define also two involutions acting on sentences: reversal and complement. The reversal of $I$, denoted by $I^r$ is obtained by writing the words of $I$ in the reverse order
\begin{equation}
(w_1,w_2,\dots , w_m)^r = (w_m, w_{m-1} , \dots , w_1) .
\end{equation}
The complement of $I$, denoted by $I^c$ is the sentence with the same maximal word as $I$ but whose image under the map $\mathrm{set}$ is the complementary subset of $\{ 1,2,\dots ,|I|-1\}$. Equivalently, in the maximal word $w(I)$ we put separating commas between letters if there was no comma between the letters in $I$.
\begin{Ex}
For $I=(aba,ca,bac)$ we have $I^r = (bac,ca,aba)$ and $I^c = (a,b,ac,ab,a,c)$
\end{Ex}
We may represent sentences in terms of ribbons placing its words in subsequent rows such that the first letter of the next word is exactly below the last letter of the previous one. Then the ribbon diagram of $I^c$ is obtained by transposition of the ribbon diagram of $I$, see Figure~\ref{fig:ribbons}.
\begin{figure}[h!]
	\begin{center}
\vskip1cm \hskip-3cm	$ \gyoung(;a;b;a,::;c;a,:::;b;a;c)$ \vskip-2.5cm \hskip5cm	
	 $ \gyoung(a,b,ac,:;ab,::;a,::;c) $	
	 \end{center}
\caption{The ribbon diagrams  of the sentence $(aba,ca,bac)$ and of its complement}
\label{fig:ribbons}	
\end{figure}

From yet another point of view, we may identify words with the so called planted tall trees (or ladders) 
	\begin{equation}
	w = a_{i_1} a_{i_2} \dots  a_{i_k} \longleftrightarrow (B_{i_k}^+ \circ \dots \circ B_{i_1}^+) (\ttree).
	\end{equation}
	Then sentences are in correspondence with concatenations of such trees, see Figure~\ref{fig:tall-tree} for an example.
\begin{Cor} \label{cor:NSym-CF}
	The set of sentences over $A$ is bijective with the context-free language~\cite{Sudkamp} generated by grammar with
	 \begin{itemize}
		\item terminal symbols $A\cup \bar{A}$, where $\bar{A}$ is the disjoint copy of $A$ with elements $\bar{a}$ for $a\in A$,
		\item nonterminal symbols $\{ X,Z \} $ with $X$ being the initial symbol of the grammar,
		\item production rules
		\begin{equation}
		X \to 1 \, | \, XZ, \qquad Z \to \bar{a}_i Z a_i \, | \, \bar{a}_i a_i , \qquad \text{where} \quad a_i\in A.
		\end{equation}
	\end{itemize}
Roughly speaking, the second rule produces words, while the first rule builds sentences from words. The relation of the Hopf algebra of trees and its subalgebras to context-free languages will be presented in another publication. 
\end{Cor}
\begin{Ex}
	The element of the context-free language described above, which corresponds to the sentence $(aba,ca,bac)$ is $\bar{a}\bar{b}\bar{a}aba\bar{a}\bar{c}ca\bar{c} \bar{a}\bar{b}bac$.
\end{Ex}

Given two sentences $I=(w_1, \dots , w_m)$, $J=(v_1,\dots ,v_n)$ their concatenation $I\cdot J$ is the sentence $(w_1, \dots , w_m,v_1,\dots ,v_n)$ obtained by juxtaposition of the sequences of their words, and corresponds to joining bunches of tall trees. 
Their near-concatenation $I\odot J$ is the sentence $(w_1, \dots , w_m v_1,\dots ,v_n)$ in which the last word of $I$ is concatenated with the first word of $J$ For example
\begin{equation}
(aba,cb) \cdot (bb,ac) = (aba,cb,bb,ac), \qquad  (aba,cb) \odot (bb,ac) = (aba,cbbb,ac).
\end{equation}
Notice that since an admissible subtree of a tall tree and its complementary tree are also of such form then the pruning coproduct \eqref{eq:Delta-cop} of tall trees doesn't lead out of that space.
\begin{Prop}
The planted tall trees generate Hopf subalgebra of ROC-trees with the concatenation multiplication and pruning coproduct.
\end{Prop}
	\begin{figure}[h!] 
	$(aba,ca,bac) \longleftrightarrow $
	\begin{picture}(50,60)(-30,0) 
	\color{red}
	\put(8,-2){\line(-1,1){20}}
	\put(-11,4){$a$}
	\put(-12,38){\line(0,1){20}}
	\put(-20,46){$a$}
	\put(8,-2){\line(0,1){20}}	
	\put(1,8){$a$}
	\put(28,18){\line(0,1){20}}
	\put(20,26){$a$}
	\color{blue}
	\put(-12,18){\line(0,1){20}}
	\put(-20,26){$b$}
	\put(28,38){\line(0,1){20}}
	\put(20,46){$b$}	
	\color{green}
	\put(8,18){\line(0,1){20}}
	\put(1,26){$c$}
	\put(8,-2){\line(1,1){20}}
	\put(20,4){$c$}	
	\color{black}
	\put(8,-2){\circle*{5}}
	\put(-12,18){\circle*{5}}
	\put(-12,38){\circle*{5}}
	\put(-12,58){\circle*{5}}
	\put(8,18){\circle*{5}}
	\put(8,38){\circle*{5}}
	\put(28,18){\circle*{5}}
	\put(28,38){\circle*{5}}
	\put(28,58){\circle*{5}}
	\end{picture} 
	\caption{Sentences as coloured tall trees}
	\label{fig:tall-tree}
\end{figure}
To describe the corresponding comultiplication in the language of sentences (coloured compositions) let us consider also weak sentences, which may contain empty words. Given weak sentence $I$, by $\widetilde{I}$ we denote the corresponding sentence obtained by removing the empty words from it (or the empty word if $I$ consists of the empty words only). We say that the weak sentence $J = (v_1, \dots , v_m)$ is contained in the sentence $I = (w_1, \dots ,w_m)$, denoted by $J\subset I$, if there exists complementary weak sentence $(u_1,\dots ,u_m)$, denoted by $I\setminus J$, such that $w_i = u_i v_i$ for $i=1,\dots ,m$. 
\begin{Rem}
	To avoid confusion we recall that the empty sentence, in accordance to our previous notation, is denoted by by $1$, the corresponding empty composition was denoted previously by $\emptyset$, and the corresponding trivial tree by $\ttree$.
\end{Rem}
\begin{Prop}
	The algebra of sentences (coloured compositions) is a bialgebra with the multiplication being concatenation of sentences, the comultiplication given by
	\begin{equation} \label{eq:cop-NS}
	\Delta(I) = \sum_{J\subset I} \widetilde{I\setminus J}\otimes \widetilde{J},
	\end{equation}
	the natural unity map, and the counit 
	\begin{equation}
	\epsilon(I)
 = \begin{cases} 1 & \text{if} \quad I = 1, \\ 0  &  \text{otherwise} \end{cases}
	\end{equation}
\end{Prop}

\begin{Ex}
	The coproduct of a ROC-tree visualized in Figure~\ref{fig:coproduct} reads in the present setting as follows 
	\begin{align}
	\Delta(ab,b) & = (ab,b)\otimes 1 + (ab) \otimes (b) + (a,b)\otimes (b) + (a) \otimes (b,b) + (b)\otimes (ab) + 1\otimes (ab,b) = \nonumber \\
	& = \left(  (ab) \otimes 1 + (a) \otimes (b) + 1 \otimes (ab)  \right) .
	 \left(  (b) \otimes 1 + 1 \otimes (b)\right).
	\end{align}
\end{Ex}
\begin{Rem}
	As the above Example shows, the comultiplication in the Hopf algebra of sentences is not cocommutative.
\end{Rem}
\begin{Cor} \label{cor:dimension}
	The bialgebra of sentences over alphabet $A$ is graded, with the weight of a sentence being its size, locally finite and connected (thus Hopf algebra). The dimension of the graded component consisting of sentences of size $m>0$ is $|A|^m 2^{m-1}$.
\end{Cor}
\begin{Cor} \label{cor:Delta-w}
	The action of the comultiplication on single-word generators reads as follows
	\begin{equation} \label{eq:Delta-w}
	\Delta(a_{i_1} a_{i_2} \dots a_{i_k}) 
	%& = 1\otimes (a_{i_1} a_{i_2} \dots a_{i_k}) +
	%(a_{i_1}) \otimes (a_{i_2} \dots a_{i_k}) + \dots + (a_{i_1} a_{i_2} \dots %a_{i_k}) \otimes 1 = \nonumber \\ &
	 = \sum_{j=0}^k (a_{i_1} a_{i_2} \dots a_{i_j}) \otimes 
	(a_{i_{j+1}}  \dots a_{i_k}) .
	\end{equation}
\end{Cor}
\begin{Rem}
	In the present setting the free Hopf algebra of Example~\ref{ex:FHA} should be identified with the subalgebra of sentences built out of single-letter words, see also Corollary~\ref{cor:kA-t}.
\end{Rem}
In \cite{Foissy-I} one can find also detailed description of the antipode of the Hopf algebra of ROD-forests, which can be used to define the antipode of the Hopf algebra of tall trees, and thus to transfer it into the language of the Hopf algebra of sentences. To make the paper self-contained we perform below the corresponding calculation from scratch avoiding this route.
\begin{Prop}
The antipode in the Hopf algebra of sentences is given by the following formula
\begin{equation}
S(I) = \sum_{ J \preccurlyeq I^r } (-1)^{\ell(J)} J .
\end{equation}	
\end{Prop}
\begin{proof}
	We will show first that the above formula gives the antipode for single-word sentences, which generate the algebra of sentences. The coproduct formula~\eqref{eq:Delta-w} and equation \eqref{eq:antipode-S} give the recurrence relation 
	\begin{equation} \label{eq:S-Hw}
	S(a_{i_1} a_{i_2} \dots  a_{i_{k}})= - \sum_{j=0}^{k-1} 
	S(a_{i_1} a_{i_2} \dots  a_{i_j}) \cdot (a_{i_{j+1}}  \dots  a_{i_{k}}),
	\end{equation}
which, in particular, for $k=1$ gives the correct formula
	\begin{equation}
	S(a_i) = - (a_i) .
	\end{equation}
	
	Assume that the expression for the antipode holds true for generators indexed by single-word sentences of size not greater than $k$, then for $k+1$ we have 
	\begin{equation}
	S(a_{i_1} a_{i_2} \dots  a_{i_{k+1}}) = 
	% &- \sum_{j=0}^k S(H_{(a_{i_1} a_{i_2} \dots  a_{i_j})} ) 
	%H_{(a_{i_{j+1}}  \dots  a_{i_{k+1}})}
	\sum_{j=0}^k \left( \sum_{J \preccurlyeq (a_{i_1} a_{i_2} \dots  a_{i_j})}(-1)^{\ell(J) + 1} J \cdot (a_{i_{j+1}}  \dots  a_{i_{k+1}})\right).
	\end{equation}
	which gives the correct expression, because we separated the last word of the sentence refining $(a_{i_1} a_{i_2} \dots  a_{i_{k+1}})$.
	
	By the anti-endomorphism property of the antipode we have
	\begin{align} \nonumber
	S(w_1,\dots , w_m) & = S(w_m) \cdot \ldots \cdot S(w_1)
	= \sum_{J_m\preccurlyeq (w_m), \ldots J_1 \preccurlyeq (w_1)} (-1)^{\ell(J_m) + \dots + \ell(J_1)} J_m \cdot \ldots \cdot J_1 = \\
	&= \sum_{J_m\dots J_1 \preccurlyeq (w_m,\dots , w_1)} (-1)^{\ell(J_m) + \dots + \ell(J_1)} J_m \cdot \ldots \cdot J_1 = \sum_{J \preccurlyeq (w_1,\dots , w_m)^r} (-1)^{\ell(J)} J,
	\end{align}
	what concludes the proof.
\end{proof}

\begin{Rem}
In	\cite{BaumannHohlweg,BergeronHohlweg,HsiaoPetersen,MantaciReutenauer,Poirier} another notion of \emph{coloured compositions} is considered. In our approach such a variant corresponds to sentences made of words with definite colours, for example $(bbb,a,bb,cccc)$. Because concatenation and pruning operations leave such property untouched one obtains this way a Hopf subalgebra of that introduced in this Section. Notice~\cite{HsiaoPetersen} that the dimension of the graded component consisting of such compositions/sentences of size $m>0$ is $|A|(|A|+1)^{m-1}$, to be compared with the dimension calculated in Corollary~\ref{cor:dimension}. 
\end{Rem}
\subsection{Coloured non-commutative symmetric functions} 
Because of the isomorphism of the Hopf algebra of sentences on unary alphabet with the Hopf algebra $\mathrm{NSym}$ of non-commutative symmetric functions, the algebra of sentences over alphabet $A$ can be also called the algebra of coloured non-commutative symmetric functions, and denoted by $\mathrm{NSym}_A$. We will discuss also other bases of $\mathrm{NSym}_A$ indexed by sentences, therefore the linear basis of sentences will be denoted from now on by $(H_I)$ and called the basis of complete homogeneous coloured non-commutative symmetric functions. The multiplication, comultiplication and the antipode in the new notation read
\begin{align}
H_I \cdot H_J & = H_{I\cdot J}, \label{eq:HI-HJ}\\ 
\Delta(H_I) & = \sum_{J \subset I} H_{\widetilde{I\setminus J}} \otimes 
H_{\widetilde{J}} ,\label{eq:Delta-HI} \\
S(H_{I}) & = \sum_{ J \preccurlyeq I^r } (-1)^{\ell(J)} H_J . \label{eq:S-HI}
\end{align}
In particular, for functions indexed by single-word sentences we have 
\begin{equation} \label{eq:Delta-H}
\Delta(H_{(a_{i_1} a_{i_2} \dots  a_{i_k})}) = \sum_{j=0}^k 
H_{(a_{i_1} a_{i_2} \dots  a_{i_j})} \otimes H_{(a_{i_{j+1}}  \dots  a_{i_k})} \;,
\end{equation}
\begin{equation} \label{eq:S-Hw}
S(H_{(a_{i_1} a_{i_2} \dots  a_{i_{k}})}) = - \sum_{j=0}^{k-1} 
S(H_{(a_{i_1} a_{i_2} \dots  a_{i_j})} ) \cdot H_{(a_{i_{j+1}}  \dots  a_{i_{k}})} .
\end{equation}
\begin{Rem}
	We define the operations of reversal and complement in the basis $(H_I)$ 
	\begin{equation} \label{eq:H-rc}
	r(H_I) = H_{I^r}, \qquad c(H_I) = H_{I^c},
	\end{equation}
	and extend them to $\mathrm{NSym}_A$ by linearity.
\end{Rem}
\begin{Ex}
	\begin{equation}
S(H_{(ab,c)}) = S(H_{(c)}) \cdot S(H_{(ab)}) = (- H_{(c)})\cdot ( -H_{(ab)} + H_{(a,b)}) = H_{(c,ab)} - H_{(c,a,b)}.
	\end{equation}	
\end{Ex}

Like in the classical case define coloured non-commutative elementary symmetric functions by
\begin{equation} \label{eq:E-def}
E_I = \sum_{J \preccurlyeq I} (-1)^{|I| - \ell(J)} H_J,
\end{equation}
what allows to rewrite the antipode of the complete homogeneous functions as
\begin{equation} \label{eq:S-H}
S(H_I) = (-1)^{|I|}E_{I^r}.
\end{equation}
\begin{Rem}
	In general $E_{I^r} \neq r(E_I)$, for example
	\begin{equation*}
	E_{(ab,c)^r} = H_{(c,a,b)} - H_{(c,ab)} , \quad \text{while} \quad
	r(E_{(ab,c)}) =  H_{(c,b,a)} - H_{(c,ab)}.
	\end{equation*}
\end{Rem}
\begin{Prop} \label{prop:mult-E}
	The product of coloured non-commutative elementary symmetric functions satisfies the formula
	\begin{equation}
	E_I \cdot E_J = E_{I\cdot J},
	\end{equation}
	in particular, the elementary symmetric functions are generated by single-word elementary functions.
\end{Prop}
\begin{proof} By equation~\eqref{eq:S-H} and anti-endomorphism property of the antipode we have
	\begin{equation}
		E_I \cdot E_J = (-1)^{|I|+|J|} S(H_{I^r}) \cdot S(H_{J^r}) = (-1)^{|I\cdot J|} S(H_{J^r \cdot I^r}) =  (-1)^{|I \cdot J|} S(H_{(I\cdot J)^r}) = E_{I \cdot J}.
	\end{equation}
\end{proof}
\begin{Prop}
	The coloured non-commutative elementary symmetric functions form a linear basis of the Hopf algebra $\mathrm{NSym}_A$, in particular
	\begin{equation} \label{eq:H-E}
	H_I = \sum_{J \preccurlyeq I} (-1)^{|I| - \ell(J)} E_J .
	\end{equation}
\end{Prop}
\begin{proof} The right hand side of equation \eqref{eq:H-E} reads
\begin{equation}
%\sum_{J \preccurlyeq I} (-1)^{|I| - \ell(J)} E_J = 
\sum_{J \preccurlyeq I} (-1)^{|I| - \ell(J)} \sum_{K \preccurlyeq J} (-1)^{|J| - \ell(K)} H_K \stackrel{|I|=|J|}{=}  \sum_{K \preccurlyeq I} (-1)^{\ell(I) - \ell(K)} H_K \sum_{K \preccurlyeq J \preccurlyeq I} (-1)^{\ell(I) - \ell(J)}.
\end{equation}	
By properties of the Moebius function \eqref{eq:Moebius} the second sum equals $1$ for $K=I$ and vanishes otherwise, what concludes the proof.
\end{proof}
\begin{Cor}
The single-word coloured non-commutative elementary symmetric functions generate the algebra $\mathrm{NSym}_A$.
\end{Cor}

Because for $|A|>1$ the Hopf algebra $\mathrm{NSym}_A$ is both non-commutative and non-cocommutative we cannot expect that the antipode is an involution. It turns out that its superposition with the reversal is.
\begin{Prop}
	In the Hopf algebra $\mathrm{NSym}_A$ of coloured non-commutative symmetric functions 
	\begin{equation} \label{eq:SrSr}
	S\circ r \circ S \circ r = \mathrm{id}.
	\end{equation}
\end{Prop}
\begin{proof}
		Notice first that
		\begin{equation} \label{eq:H-SrE}
		(S\circ r)(E_I)  = S\left(   \sum_{J \preccurlyeq I} (-1)^{|I| - \ell(J)} H_{J^r}       \right) =  (-1)^{|I|} \sum_{J \preccurlyeq I} (-1)^{ |J_r| - \ell(J)} E_{J} =  (-1)^{|I|} H_I,
		\end{equation}
		where we used equations~\eqref{eq:S-H}, \eqref{eq:H-E} and the fact that $|J^r|=|I|$ for $J \preccurlyeq I$.
Then in the basis of complete symmetric functions by equations \eqref{eq:S-H} and \eqref{eq:H-SrE} we have
	\begin{equation}
	(S\circ r \circ S \circ r )(H_I) = (-1)^{|I|} (S\circ r)(E_I) = H_I,
	\end{equation}
what concludes the proof.
\end{proof}

We will close this Section by presenting coloured and non-commutative analogs of some classical properties of single-word elementary symmetric functions. Notice that equations \eqref{eq:antipode-S} and \eqref{eq:E-def} can be rewritten for $k>0$ as
\begin{equation}
\sum_{j=0}^k 
(-1)^j E_{(a_{i_1} a_{i_2} \dots  a_{i_j})}  H_{(a_{i_{j+1}}  \dots  a_{i_k})} =
\sum_{j=0}^k 
(-1)^j H_{(a_{i_1} a_{i_2} \dots  a_{i_j})}  E_{(a_{i_{j+1}}  \dots  a_{i_k})} = 0.
\end{equation}
In particular, we have two recurrence formulas which start from $E_1 = 1$, and read
\begin{align} \label{eq:E-rec}
E_{(a_{i_1} a_{i_2} \dots  a_{i_k})} = \sum_{j=0}^{k-1} 
(-1)^{k-j+1} E_{(a_{i_1} a_{i_2} \dots  a_{i_j})}  H_{(a_{i_{j+1}}  \dots  a_{i_k})} =
\sum_{j=1}^k 
(-1)^{j-1} H_{(a_{i_1} a_{i_2} \dots  a_{i_j})}  E_{(a_{i_{j+1}}  \dots  a_{i_k})} .
\end{align}
Finally we present formula for their coproduct.
\begin{Prop}
	The analog of the coproduct formula~\eqref{eq:Delta-H} but for coloured non-commutative elementary symmetric functions reads as follows
\begin{equation}
\Delta(E_I) = \sum_{J \subset I} E_{\widetilde{J}} \otimes 
E_{\widetilde{I\setminus J}} ,
\end{equation}
which in the case of single-word functions gives
\begin{equation} \label{eq:Delta_Ew}
\Delta( E_{(a_{i_1} a_{i_2} \dots  a_{i_k})}  ) = \sum_{j=0}^k 
E_{(a_{i_{j+1}}  \dots  a_{i_k})} \otimes  E_{(a_{i_1} a_{i_2} \dots  a_{i_j})} .
\end{equation}
\end{Prop}
\begin{proof}
It is enough to prove the second equation, because then Proposition~\ref{prop:mult-E} implies the first one. Since it holds for $k=1$, then we can start induction by applying the coproduct operation on the recurrence \eqref{eq:E-rec}. Using the homomorphism property of the comultiplication we can expand corresponding expressions and collect coefficients at various terms of consecutive degrees on the right hand side of the tensor product sign. By the recurrence relations  \eqref{eq:E-rec} most of them vanishes, and what remains gives equation~\eqref{eq:Delta_Ew}.
\end{proof}
\begin{Ex}
	To calculate the coproduct of $E_{(ab)} = E_{(a)}H_{(b)} - H_{(ab)}$ first notice that
	\begin{equation*}
	\Delta(E_{(ab)}) = (E_{(a)}\otimes 1 + 1 \otimes E_{(a)})(H_{(b)}\otimes 1 + 1 \otimes H_{(b)}) - ( H_{(ab)} \otimes 1 + H_{(a)} \otimes H_{(b)} + 1 \otimes H_{(ab)}).
	\end{equation*}
	There is only one term $(\dots)\otimes 1$ of the right degree zero. Its coefficient is $E_{(a)} H_{(b)} - H_{(ab)} = E_{(ab)}$. 
	There are two terms $(\dots)\otimes E_{(a)}$ and $(\dots)\otimes H_{(b)}$ of the right degree one. The coefficient of the first one equals $H_{(b)} = E_{(b)}$,
	while the coefficient of the second term reads $E_{(a)} - H_{(a)} = 0$.
Finally, the coefficients of the right degree two have on the left $1\otimes (\dots)$ and sum up to $E_{(a)}H_{(b)} - H_{(ab)} = E_{(ab)} $.
\end{Ex}
\begin{Rem}
	Notice that, contrary to the unary (monochromatic) case $|A|=1$, for $|A|>1$ the coproduct formulas for single-word coloured non-commutative complete and elementary symmetric functions are not the same.
\end{Rem}

\section{Coloured quasi-symmetric functions} \label{sec:CQSym}
In this Section we study basic properties of the graded dual to the Hopf algebra $\mathrm{NSym}_A$, which we later will call the Hopf algebra of coloured quasi-symmetric functions, and denote by $\mathrm{QSym}_A$. In particular, we introduce the dual  basis to complete function basis $(H_I)$, which will be called later the basis of coloured monomial quasi-symmetric functions. Then we will provide a realization of the algebra $\mathrm{QSym}_A$ in terms of series of bounded degree with partially commuting variables. 
\subsection{The graded dual of $\mathrm{NSym}_A$}
In the the graded dual $(\mathrm{NSym}_A)^{gr}$ of the Hopf algebra of $A$-coloured non-commutative symmetric functions, by $(H^*_I)$ denote the dual basis to the basis 
$(H_J)$ of complete symmetric functions 
\begin{equation}
\langle H^*_I , H_J \rangle = \delta_{I,J}.
\end{equation}
The dual to the concatenation product \eqref{eq:Delta-HI} is the deconcatenation coproduct $\delta$, which is given by
\begin{equation} \label{eq:delta-MI} \begin{split}
\delta (H^*_I) & = \sum_{ I = J\cdot K} 
H^*_J \otimes H^*_K \quad \text{i.e.}\\
\delta (H^*_{(w_1, w_2, \dots w_k)}) & = \sum_{j=0}^k H^*_{(w_1, w_2, \dots w_j)} \otimes 
H^*_{(w_{j+1},  \dots w_k)} \; ,
\end{split} \end{equation}
In particular, elements of the dual basis indexed by single-word sentences are primitive  elements of the coproduct.

The product in $(\mathrm{NSym}_A)^{gr}$ can be defined directly by dualization of the pruning coproduct of tall trees described in the basis of complete functions by formula~\eqref{eq:Delta-HI}. Equivalently, it can be described in terms of the original grafting product of trees and the dual to the injection map $\mathrm{NSym}_A \hookrightarrow \Bbbk T_{A}$. It is therefore restriction of the grafting product from ROC-trees to tall trees, i.e. we can graft planted tall trees at the root or on the top of another such tree, with the restriction that two trees cannot be grafted on the same top. 

\begin{figure}
\begin{center}
	\includegraphics[width=2cm]{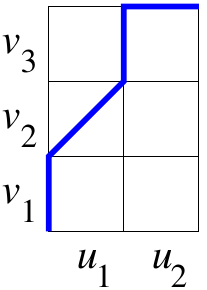}
\end{center}	
	\caption{The lattice for quasi-shuffle product of $(u_1,u_2)$ and $(v_1,v_2,v_3)$ with path representing the summand $(v_1,u_1 v_2, v_3, u_2)$}
	\label{fig:Q-shuffle-lattice}
\end{figure}

The geometric procedure on the level of trees is the same as that in the monochromatic case, so we keep the same name and symbol of the quasi-shuffle product. Quasi-shuffle 
$I \Qshuffle J$ of two sentences $I = (u_1, u_2, \dots , u_k)$ and $J = (v_1, v_2, \dots , v_m)$ is thus the sum of shuffles of components $u_i$ and $v_j$ of $I$ and $J$, where in addition we may replace any number of pairs of consecutive words $u_i, v_j$ in the shuffle by their concatenation $u_i v_j$. It can be represented, compare with~\cite{Introduction-QSym}, by a path in the lattice of size $k\times m$ from the left bottom corner to the right top corner with horizontal steps $(1,0)$ representing words $u_i$, vertical steps $(0,1)$ representing words $v_j$ and oblique steps $(1,1)$ representing words $u_i v_j$, see Figure~\ref{fig:Q-shuffle-lattice} for an example.

Therefore we have (compare with equation~\eqref{eq:M-a-M-b})
\begin{equation} \label{eq:M-I-M-J}
H^*_I H^*_J = \sum_{K} H^*_K \;,
\end{equation}
where $K$ can be obtained as a summand in quasi-shuffle of $I$ and $J$. By dualizing the coproduct formula~\eqref{eq:Delta-HI} we can see that there exists week sentence $J^\prime \subset K$ such that $J=\widetilde{J^\prime}$ and $I = \widetilde{K\setminus J^\prime}$.
\begin{Rem}
	Similarly to the quasi-shuffle product of compositions, the quasi-shuffle of sentences can be defined recursively for all sentences $I,J$ and all non-empty words $u,v$ by
	\begin{equation*}
	I\Qshuffle 1 = 1 \Qshuffle I = I, \quad ((u)\cdot I) \Qshuffle ((v)\cdot J) = (u)\cdot(I\Qshuffle ((v)\cdot J)) + (v)\cdot(((u)\cdot I)\Qshuffle J)+ (uv)\cdot (I\Qshuffle J).
	\end{equation*}
\end{Rem}
\begin{Ex} \label{ex:q-sh-H*}
The quasi-shuffle product of two tall trees takes the form given in Figure~\ref{fig:q-shuffle-prod}, and the product of corresponding monomial functions reads
\begin{equation} \label{eq:q-shuffle-prod}
H^*_{(a,b)} H^*_{(b)} =
2 H^*_{(a,b,b)} + H^*_{(a,bb)} +
 H^*_{(ab,b)} +
H^*_{(b,a,b)} \; ,
\end{equation}
compare also with Figure \ref{fig:ashuffle-prod} describing the asymmetric shuffle product of the same trees.
\begin{figure}[h!]
\begin{picture}(30,50)(-30,0) 
\color{red}
\put(-12,0){\line(1,2){10}}
\put(-4,8){$b$}
\color{blue}
\put(-12,0){\line(-1,2){10}}
\put(-26,8){$a$}
\color{black}
\put(-12,0){\circle*{5}}
\put(-2,20){\circle*{5}}
\put(-22,20){\circle*{5}}
\end{picture}
$\Qshuffle$
\begin{picture}(20,40)
\color{red}
\put(8,0){\line(0,1){20}}
\put(12,8){$b$}
\color{black}
\put(8,0){\circle*{5}}
\put(8,20){\circle*{5}}
\end{picture}
$=$
\begin{picture}(50,40)
\color{red}
\put(27.5,0){\line(0,1){20}}
\put(28.5,0){\line(0,1){20}}
\put(31,9){$b$}
\put(28,0){\line(1,1){20}}
\put(44,7){$b$}
\color{blue}
\put(27.5,0){\line(-1,1){20}}
\put(28.5,0){\line(-1,1){20}}
\put(8,7){$a$}
\color{black}
\put(28,0){\circle*{5}}
\put(28,20){\circle*{5}}
\put(8,20){\circle*{5}}
\put(48,20){\circle*{5}}
\end{picture}
$+$
\begin{picture}(30,40)
\color{red}
\put(18,0){\line(1,2){10}}
\put(27.5,20){\line(0,1){20}}
\put(28.5,20){\line(0,1){20}}
\put(28,8){$b$}
\put(32,29){$b$}
\color{blue}
\put(17.5,0){\line(-1,2){10}}
\put(18.5,0){\line(-1,2){10}}
\put(3,8){$a$}
\color{black}
\put(8,20){\circle*{5}}
\put(28,40){\circle*{5}}
\put(18,0){\circle*{5}}
\put(28,20){\circle*{5}}
\end{picture}
$+$
\begin{picture}(50,40)
\color{red}
\put(28,0){\line(0,1){20}}
\put(31,9){$b$}
\put(27.5,0){\line(1,1){20}}
\put(28.5,0){\line(1,1){20}}
\put(44,7){$b$}
\color{blue}
\put(27.5,0){\line(-1,1){20}}
\put(28.5,0){\line(-1,1){20}}
\put(8,7){$a$}
\color{black}
\put(28,0){\circle*{5}}
\put(28,20){\circle*{5}}
\put(8,20){\circle*{5}}
\put(48,20){\circle*{5}}
\end{picture}
+
\begin{picture}(30,40)
\color{red}
\put(18,0){\line(-1,2){10}}
\put(17.5,0){\line(1,2){10}}
\put(18.5,0){\line(1,2){10}}
\put(3,8){$b$}
\put(28,8){$b$}
\color{blue}
\put(7.5,20){\line(0,1){20}}
\put(8.5,20){\line(0,1){20}}
\put(12,29){$a$}
\color{black}
\put(8,20){\circle*{5}}
\put(8,40){\circle*{5}}
\put(18,0){\circle*{5}}
\put(28,20){\circle*{5}}
\end{picture}
$+$
\begin{picture}(40,40)
\color{red}
\put(28,0){\line(-1,1){20}}
\put(27.5,0){\line(1,1){20}}
\put(28.5,0){\line(1,1){20}}
\put(8,7){$b$}
\put(44,7){$b$}
\color{blue}
\put(27.5,0){\line(0,1){20}}
\put(28.5,0){\line(0,1){20}}
\put(31,9){$a$}
\color{black}
\put(28,0){\circle*{5}}
\put(28,20){\circle*{5}}
\put(8,20){\circle*{5}}
\put(48,20){\circle*{5}}
\end{picture}
\caption{The quasi-shuffle product of two tall trees as given in equation \eqref{eq:q-shuffle-prod}; grafted planted tall trees are thickened}
\label{fig:q-shuffle-prod}
\end{figure}
\end{Ex}
\begin{Rem}
We have equipped the space of sentences over $A$ with the dual Hopf algebra structure, graded dual to that described in Section~\ref{sec:sent-A}, with quasi-shuffle product and deconcatenation coproduct. Being dual to non-commutative and non-cocommutative Hopf algebra the new algebra, for $|A|>1$, is also both non-commutative and non-cocommutative.
\end{Rem}
Let us calculate the number of quasi-shuffle paths with prescribed number of oblique steps, which will be used in Section~\ref{sec:power-sum}.

\begin{Prop} \label{prop:q-shuffle-paths-i}
	The number of quasi-shuffle paths in $k\times m$ lattice with exactly $i$ oblique steps $(1,1)$, $i=0, 1, \dots , \min \{ k,m \}$, equals
	\begin{equation} \label{eq:k-m-i}
	\binom{k + m - 2i}{k-i} \binom{k+m-i}{i} = 
	\binom{k }{i} \binom{k+m-i}{k} = \binom{m}{i} \binom{k+m-i}{m}	.
	\end{equation}	
\end{Prop}
\begin{proof}
	Two quasi-shuffle paths are called to have the same shuffle part if they coincide after contracting all the oblique segments $(1,1)$, see Figure~\ref{fig:Q-shuffle-part} for an example.
	\begin{figure}
		\begin{center}
			\includegraphics[width=8 cm]{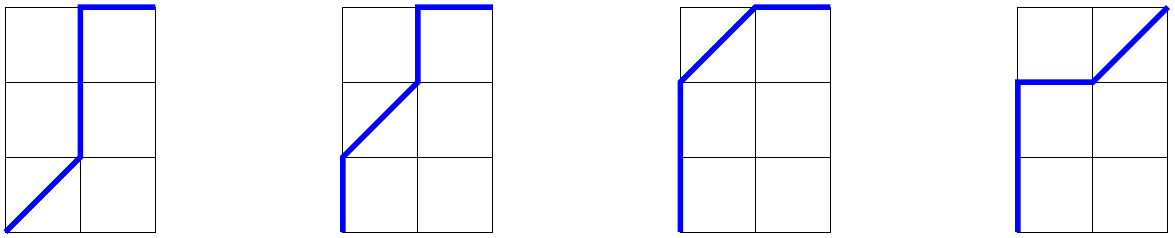}
		\end{center}	
		\caption{Four quasi-shuffle paths in $2\times 3$ lattice which have the same shuffle part. The difference is in location of the $(1,1)$ segment which can be placed at any vertex of the contracted path}
		\label{fig:Q-shuffle-part}
	\end{figure}
	Decomposition of the set of paths with exactly $i$  oblique segments into disjoint classes having the same shuffle part, and then fixing location of the segments, gives the first part of formula \eqref{eq:k-m-i}. Two other expressions, which can be derived by simple algebra, also have combinatorial interpretation .	
	The second/third one means that the path can be encoded by first fixing columns/rows for the oblique segments $(1,1)$, and then by choosing which steps of the path are vertical/horizontal segments $(0,1)$/$(1,0)$. 
\end{proof}
\subsection{Generalized quasi-symmetric functions} \label{sec:Q-Sym}
For each colour $a\in A$, let $x_a = ( x_{a,1}, x_{a,2}, x_{a,3} , \dots )$ denote infinite totally ordered set of variables, each of degree $1$, define also
$x_A = \bigcup_{a\in A}x_a$. We assume \emph{partial commutativity} of the variables, i.e. within each set the variables commute, but for different colours commutativity is allowed for different second indices only
\begin{equation} \label{eq:C-F-p-c}
x_{a,i} x_{b,j} = x_{b,j} x_{a,i} \qquad \text{for} \qquad i\neq j, \qquad a,b\in A ,
\end{equation}
%\begin{Rem}
%In quantum physics, requirement of the form "elements with different indices commute" %is referred often as the \emph{ultralocality condition}.
%\end{Rem}
%\begin{Cor}
%Notice that within each set (the first index is fixed) all variables commute. In %particular, for $n=1$ we obtain totally ordered set of commuting variables.
%\end{Cor}

We will consider a subset $\mathrm{QSym}_A$ of the algebra $\Bbbk[x_A]$ of series of bounded degree with natural multiplication, which can be described as follows. Due to partial commutativity any monomial in variables $x_{a,i}$ can be uniquely reordered in such a way that the second indices of variables form weakly increasing (finite) sequence, say 
$ i_1\leq i_2 \leq \dots \leq i_k $. Given word $w=a_{i_1}a_{i_2}\dots a_{i_k}$ and given $j\in\mathbb{N}$, by $x_{w,j}$ denote monomial of degree $|w|$
\begin{equation}
x_{w,j} = x_{a_{i_1},j}x_{a_{i_2},j}\dots x_{a_{i_k},j},
\end{equation}
for example $x_{abb,2} = x_{a,2} x_{b,2}^2$. 
From the other side, given reordered monomial $x_{w_1,j_1} x_{w_2,j_2} \dots x_{w_m,j_m}$ with $j_1 < j_2 < \dots < j_m$ by its sentence we mean $(w_1,w_2,\dots ,w_m)$. A formal series belongs to $\mathrm{QSym}_A$ when its coefficients in front of monomials with the same sentence coincide.

It is easy to see that the set $\mathrm{QSym}_A$ is in fact linear space with basis indexed by sentences. In fact, given sentence $I=(w_1,w_2,\dots ,w_m)$, by $M_I$ denote the infinite series of the finite degree $|I|$
\begin{equation}
M_I = \sum_{1\leq j_1 < j_2 < \dots <j_m} x_{w_1,j_1} x_{w_2,j_2} \dots x_{w_m,j_m},
\end{equation}
which will be called a coloured monomial quasi-symmetric function. 
\begin{Ex}
	Consider product of two such series 
	\begin{equation*}
	M_{(a,b)}  = x_{a,1} x_{b,2} + x_{a,1} x_{b,3} + 
	x_{a,2} x_{b,3 } + \dots \qquad \text{and} \qquad 
M_{(b)} = x_{b,1} + x_{b,2} + x_{b,3} + \dots ,
	\end{equation*}
	which after the reordering reads
	\begin{equation*}
	\begin{split}
	& M_{(a,b)}  M_{(b)} =
	( x_{a,1} x_{b,2} + x_{a,1} x_{b,3} + x_{a,2} x_{b,3 } + \dots )
	( x_{b,1} + x_{b,2} + x_{2,3} + \dots ) =  \\
	= & (x_{a,1} x_{b,1} x_{b,2} + \dots ) + 
	(x_{a,1} x_{b,2} x_{b,2} + \dots ) +
	2 (x_{a,1} x_{b,2} x_{b,3} + \dots ) + 
	(x_{b,1} x_{a,2} x_{b,3} + \dots ) = \\
	= & M_{(ab,b)} + M_{(a,bb)}  + 
	2 M_{(a,b,b)}  + M_{(b,a,b)} , 
	\end{split}
	\end{equation*}
	and compare with Example \ref{ex:q-sh-H*}
	or Figure~\ref{fig:q-shuffle-prod}.
\end{Ex}
\begin{Prop} \label{prop:H-gr-Q}
	The subspace $\mathrm{QSym}_A$ of $\Bbbk[x_A]$ spanned by the series $(M_I)$ is a subalgebra isomorphic to $(\mathrm{NSym}_A)^{gr}$ with the isomorphism given by $M_I\leftrightarrow H^*_I$.
\end{Prop}
\begin{proof}
Multiplication of two monomials with different second indices gives monomials with shuffled words. When second indices of two words of both monomials coincide then in the multiplication and reordering procedure the words will be concatenated.	
\end{proof}

To define the coproduct $\delta$ in $\mathrm{QSym}_A$ we use the doubling variables trick described in Section~\ref{sec:qs-ns}, i.e. to the set of variables $x_A$ we add its copy $y_A$. We place the new variables \emph{after} the old ones, what in particular implies that the new and old variables commute and allows to separate the variables in the reordering process. 
\begin{Ex}
Applying the doubling variable procedure to $M_{(ab,b)} = x_{a,1} x_{b,1} x_{b,2}  +  \dots $ we obtain
\begin{equation*} \begin{split}
x_{a,1} x_{b,1} x_{b,2}  +  \dots & \mapsto 
x_{a,1} x_{b,1} x_{b,2} + \dots +  x_{a,1} x_{b,1} y_{b,1} +  \dots + 
y_{a,1} y_{b,1} y_{b,2} + \dots = \\
& = M_{(ab,b)}(x)  +  M_{(ab)}(x) 
M_{(b)}(y)  +
M_{(ab,b)}(y) ,
\end{split}
\end{equation*}
getting this way 
\begin{equation*}
\delta(M_{(ab,b)}) = M_{(ab,b)} \otimes 1 +
M_{(ab)} \otimes M_{(b)}  +
 1 \otimes M_{(ab,b)} \; .
\end{equation*}
\end{Ex}
\begin{Prop}
The algebra isomorphism described in Proposition \ref{prop:H-gr-Q} is the Hopf algebra isomorphism. In particular the deconcatenation coproduct in $(\mathrm{NSym}_A)^{gr}$ can be realized by the variables doubling method in $\mathrm{QSym}_A$
\begin{equation} \label{eq:copr-q-sym-M}
\begin{split}
\delta (M_I) & = \sum_{ I = J.K} 
M_J \otimes M_K \quad \text{i.e.}\\
\delta (M_{(w_1, w_2, \dots w_k)}) & = \sum_{j=0}^k M_{(w_1, w_2, \dots w_j)} \otimes 
M_{(w_{j+1},  \dots w_k)} \; .
\end{split}
\end{equation}
\end{Prop}
\begin{proof}
It is enough to consider how the variables doubling works for the monomial
$x_{w_1,1} x_{w_2,2} \dots x_{w_m,m}$. Finiteness of the size of the sentence assures finite sum decomposition.
\end{proof}
\begin{Rem}
There is no need to check coassociativity of the coproduct in $\mathrm{QSym}_A$ or its compatibility with the (quasi-shuffle) product, because this holds by definition of $(\mathrm{NSym}_A)^{gr}$.
\end{Rem}

Finally let us present the formula for the antipode $S^*$ in the monomial basis of $\mathrm{QSym}_A$. By dualizing the corresponding formula~\eqref{eq:S-HI} for the antipode $S$ in the complete basis of $\mathrm{NSym}_A$ we get
\begin{equation} \label{eq:S*-MI}
S^*(M_I) = (-1)^{\ell(I)} \sum_{J^r \succcurlyeq I} M_{J} 
\end{equation}
which for unary words (i.e. usual compositions) reduce to that found in~\cite{MalvenutoReutenauer,Ehrenborg}.
\begin{Ex}
	By the equation~\eqref{eq:S*-MI}
	\begin{equation}
	S^*(M_{(ab,c)}) = M_{(c,ab)} + M_{(abc)}.
	\end{equation}
	By the coproduct formula
	\begin{equation}
	\delta(M_{(ab,c)}) = 1\otimes M_{(ab,c)} + M_{(ab)}\otimes M_{(c)} +  M_{(ab,c)} \otimes 1 
	\end{equation}
	and the recurrence \eqref{eq:antipode-S} it should be equal to
	\begin{equation}
	- M_{(ab,c)} - S^*( M_{(ab)} ) M_{(c)} = 
	- M_{(ab,c)} + M_{(ab)} M_{(c)} = - M_{(ab,c)} + M_{(ab,c)} +
	M_{(c,ab)} + M_{(abc)}.
	\end{equation}
\end{Ex}
\begin{Rem}
	Notice that, contrary to the unary case~\cite{Grinberg-Reiner}, we cannot sum up in~\eqref{eq:S*-MI} with respect to $J \succcurlyeq I^r$. Such a summation in the above Exercise would give incorrect result $M_{(c,ab)} + M_{(cab)}$.
\end{Rem}
By dualizing equations \eqref{eq:H-rc} define the operations of reversal and complement in the monomial basis~$(M_I)$ 
\begin{equation}
r(M_I) = M_{I^r}, \qquad c(M_I) = M_{I^c},
\end{equation}
and extend them to $\mathrm{QSym}_A$ by linearity.
\section{The fundamental basis of $\mathrm{QSym}_A$ and its dual basis in $\mathrm{NSym}_A$} \label{sec:F-R}
In this Section we define and study another basis of the Hopf algebra of coloured quasi-symmetric functions $\mathrm{QSym}_A$, which is the analog of the fundamental basis~\cite{Stanley-q-sym,Ehrenborg}. Then we consider the dual basis in $\mathrm{NSym}_A$ to the coloured fundamental functions, which can be called the basis of coloured ribbon non-commutative Schur functions.

\subsection{Fundamental coloured quasi-symmetric functions}
Define the fundamental coloured quasi-symmetric functions indexed by sentences as
\begin{equation} \label{eq:def-FI}
F_I = \sum_{J \preccurlyeq I} M_J .
\end{equation}
By properties of the Moebius function~\eqref{eq:Moebius} of the poset of sentences one can invert the above relation
\begin{equation}
M_I = \sum_{J \preccurlyeq I} (-1)^{\ell(J) - \ell(I)} F_J ,
\end{equation}
what shows that the fundamental functions form a linear basis in $\mathrm{QSym}_A$.

Let us find expressions for coproduct, product and the antipode of $\mathrm{QSym}_A$ in the fundamental basis. 
\begin{Prop} \label{prop:delta-F} 
	The coproduct in the fundamental basis
	\begin{equation}
	\delta(F_I) = \sum_{\substack{ K^\prime \cdot K^{\prime\prime} = I \\
	K^\prime \odot K^{\prime\prime} = I} } F_{K^\prime} \otimes F_{K^{\prime\prime}}
	\end{equation}
	where the summation is over pairs of sentences which give the indexing sentence $I$ by concatenation or the near-concatenation.
\end{Prop} 
\begin{proof}
	By definition~\eqref{eq:def-FI} and coproduct formula~\eqref{eq:copr-q-sym-M} in the monomial basis, and grouping terms
	\begin{equation}
	\delta(F_I) = \sum_{ J\preccurlyeq I  } \left( 
	\sum_{J^\prime . J^{\prime\prime} =J }  M_{J^\prime} \otimes M_{J^{\prime\prime}} 
	\right) = \sum_{(K^\prime , K^{\prime\prime}  )}  F_{K^\prime} \otimes  F_{K^{\prime\prime}} ,
	\end{equation}
	where the sum is over the pairs $(K^\prime, K^{\prime\prime})$ which give splitting of $I$ into two parts. The segmentation may be either between words of $I$ or in the middle of a word. The first case gives 
	$I= K^\prime \cdot K^{\prime\prime} $, while the second one gives $I= K^\prime \odot K^{\prime\prime} $. 
\end{proof}
\begin{Ex} The deconcatenation coproduct of the fundamental function $F_{(ab,c)} = M_{(ab,c)} + M_{(a,b,c)}$ is given by
	\begin{equation}
	\delta ( F_{(ab,c)} ) = 1 \otimes F_{(ab,c)} + F_{(a)} \otimes F_{(b,c)} +
	F_{(ab)} \otimes F_{(c)} + F_{(ab,c)} \otimes 1.
	\end{equation}
\end{Ex}
In order to describe multiplication in the fundamental basis notice that in multiplying $F_I$ and $F_J$ we multiply $M_{I^\prime}$ and $M_{J\prime}$ for any $I^\prime \preccurlyeq I$ and $J^\prime \preccurlyeq J$. Then we group monomial functions into the fundamental ones. This procedure leads to definition of the fundamental shuffle $I\Fshuffle J$ of sentences described as follows.
\begin{enumerate}
	\item perform ordinary shuffle $\shuffle$ of letters of maximal words of both sentences,
	\item concatenate neighboring letters of words of $I$, and concatenate neighboring letters of words of~$J$,
	\item concatenate pairs of neighboring subwords of words of $I$ and $J$ (in this order),
\end{enumerate}  
which gives directly the desired formula.
\begin{Prop} \label{prop:FI.FJ}	
	Multiplication of two fundamental functions is given by
	\begin{equation}
	F_I \cdot F_J = \sum F_K,
	\end{equation}
where sentence $K$ is a summand of $I\Fshuffle J$.
\end{Prop}
\begin{Ex} To find $(ab)\Fshuffle (c,d)$ at each step of the procedure we obtain
	\begin{enumerate}
		\item $(a,b,c,d) + (a,c,b,d) + (a,c,d,b) + (c,a,b,d) + (c,a,d,b) + (c,d,a,b)$
		\item $(ab,c,d) + (a,c,b,d) + (a,c,d,b) + (c,ab,d) + (c,a,d,b) + (c,d,ab)$
		\item $(abc,d) + (ac,bd) + (ac,d,b) + (c,abd) + (c,ad,b) + (c,d,ab)$,
	\end{enumerate}
what gives
	\begin{equation}
	F_{(ab)} \cdot F_{(c,d)} = F_{(abc,d)} + F_{(ac,bd)} + F_{(ac,d,b)} + F_{(c,abd)} + F_{(c,ad,b)} + F_{(c,d,ab)} .
	\end{equation}
\end{Ex}
\begin{Rem}
	In passing to the unary alphabet we obtain the corresponding multiplication formula for the fundamental quasi-symmetric functions~\cite{Stanley,Grinberg-Reiner}. 
		Because the same structure of posets of compositions and of sentences the proof presented there can be transferred also to our context. One has to label letters of the maximal words of the two sentences by integers, whose descent sets model the separation of letters into words.
\end{Rem}
\begin{Prop} \label{prop:S*F}
	The antipode in the fundamental basis is given by 
	\begin{equation} \label{eq:S*-FI}
	S^*(F_I) = (-1)^{|I|} r(F_{I^c}).
	\end{equation}
\end{Prop}
\begin{proof}
	Expanding the fundamental function in the monomial basis and using of~\eqref{eq:S*-MI} we obtain
	\begin{equation}
	S^*(F_I) = \sum_{J \preccurlyeq I} S^*(M_J) = \sum_{J \preccurlyeq I} (-1)^{\ell(J)}  \sum_{K \succcurlyeq J} M_{K^r} = \sum_{K}  M_{K^r}  
	\sum_{ \substack{J\preccurlyeq K \\ J\preccurlyeq I}}(-1)^{\ell(J)},
	\end{equation}
	where the last sum is over sentences $J$ which refine simultaneously $K$ and $I$. By properties of the Moebius function of the poset of such refinements the sum doesn't vanish only if the poset consists of one element only, which must be therefore the sentence of one-letter words. In this case $\ell(J) = |I|$ and $K$ must be a refinement of the complement of $I$, i.e. $K\preccurlyeq I^c$. This gives
	the corresponding version of monochromatic formulas of~\cite{Ehrenborg,MalvenutoReutenauer}
	\begin{equation}
	S^*(F_I) = \sum_{K \preccurlyeq I^c} (-1)^{|I|} r(M_{K})  = (-1)^{|I|} r(F_{I^c}).
	\end{equation}
\end{proof}
\begin{Ex}
	In order to directly calculate the antipode of $F_{(ab,c)} = M_{(ab,c)}+M_{(a,b,c)}$ we first find from equation~\eqref{eq:S*-MI}
	\begin{equation}
	S^*(M_{(ab,c)}) = M_{(c,ab)} + M_{(abc)}, \qquad S^*(M_{(a,b,c)}) = - \left( M_{(c,b,a)} + M_{(c,ab)} + M_{(bc,a)} + M_{(abc)} \right),
	\end{equation}
	which summed up give $S^*(F_{(ab,c)}) = - \left( M_{(c,b,a)} +  M_{(bc,a)} \right)
	= - r(F_{(a,bc)})$ in agreement with~\eqref{eq:S*-FI}.
\end{Ex}
\begin{Rem}
	Notice that in general $r(F_I) \neq F_{I^r}$, for example 
	\begin{equation*}
	r(F_{(a,bc)}) = M_{(bc,a)} + M_{(c,b,a)}, \qquad \text{while} \quad
	F_{(a,bc)^r} = M_{(bc,a)} + M_{(b,c,a)}.
	\end{equation*}
\end{Rem}
\begin{Cor}
	By involutivity of reversal and complement operations we directly obtain the dual counterpart of formula~\eqref{eq:SrSr}
	\begin{equation}
	r\circ S^* \circ r \circ S^* = \mathrm{id}.
	\end{equation}
\end{Cor}
As an exercise we recommend for the interested Reader to perform the calculation in the monomial basis.
\subsection{Coloured ribbon non-commutative Schur functions}
Consider the basis $(R_I)$ in $\mathrm{NSym}_A$ dual to the fundamental basis $(F_I)$ in $\mathrm{QSym}_A$
\begin{equation}
\langle F_I, R_J \rangle = \delta_{I,J}.
\end{equation}
In the monochromatic case such basis was introduced in~\cite{Gelfand-NSym} as a non-commutative analog of the ribbon Schur functions~\cite{MacMahon}. 
\begin{Prop}
	The relation between the complete basis $(H_I)$ and the ribbon basis $(R_I)$ is given by
	\begin{align}
	H_I & = \sum_{J\succcurlyeq I} R_J,\\
	R_I & = \sum_{J\succcurlyeq I} (-1)^{\ell(J) - \ell(I)} H_J. \label{eq:R-H}
	\end{align}
\end{Prop}
\begin{proof}
For the first assertion, note that
\begin{equation}
H_I = \sum_{J} \langle L_J, H_I \rangle R_J = \sum_{J} \sum_{K \preccurlyeq J} \langle  M_K, H_I \rangle R_J = \sum_{J\succcurlyeq I} R_J.
\end{equation}
The second assertion follows from the first one by inclusion-exclusion.
\end{proof}
\begin{Cor}
	The dual version of Proposition~\ref{prop:delta-F} gives the multiplication formula in the ribbon basis
	\begin{equation}
	R_I R_J =  R_{I \cdot J} + R_{I\odot J}.
	\end{equation}
\end{Cor}
\begin{Cor}
	The dual version of Proposition~\ref{prop:FI.FJ} gives the comultiplication in the ribbon basis
	\begin{equation}
	\Delta(R_I)=  \sum_{J,K} \langle F_{J\sFshuffle K} , R_I \rangle R_J \otimes R_K.
	\end{equation}
\end{Cor}
\begin{Ex}
	Both calculations using the complete basis expansion~\eqref{eq:R-H} with the corresponding coproduct formula~\eqref{eq:Delta-HI} or the above Corollary and definition of the fundamental shuffle give
	\begin{equation*}
	\Delta(R_{(ab,c)}) = 1 \otimes R_{(ab,c)} + R_{(a)} \otimes R_{(b,c)} + R_{(c)}\otimes R_{(ab)} + R_{(a,c)}\otimes R_{(b)} + R_{(ac)}\otimes R_{(b)}+ R_{(ab,c)} \otimes 1.
	\end{equation*}
	In particular, we have
	\begin{align*}
	(a)\Fshuffle (b,c) & = \boldsymbol{(ab,c)} + (b,ac) + (b,c,a) , \\
	(c) \Fshuffle (ab) & = (cab) + (a,cb) + \boldsymbol{(ab,c)} , \\
	(a,c) \Fshuffle (b) & = (a,cb) + \boldsymbol{(ab,c)} + (b,a,c) , \\
	(ac)\Fshuffle (b) & = (acb) + \boldsymbol{(ab,c)} + (b,ac) .
	\end{align*}
\end{Ex}
In the monochromatic case there exists~\cite{Grinberg-Reiner} a convenient formula, of the form~\eqref{eq:S*-FI}, expressing the coproduct in the ribbon basis. Because in the coloured case taking refinements does not commute with reversal we can provide only the following result.
\begin{Cor}
	The dual version of equation~\eqref{eq:S*-FI} reads as follows
	\begin{equation}
	S\circ r (R_I) = (-1)^{|I|} R_{I^c} .
	\end{equation}
\end{Cor}
\begin{proof}
	This can be shown by dualization of equation~\eqref{eq:S*-FI}.  Let us provide also direct proof. By expressing the ribbon basis in the complete one, and by using the corresponding formula for the antipode we obtain
\begin{align*}
S\circ r (R_I) & = \sum_{J \succcurlyeq I}(-1)^{\ell(J) - \ell(I)} S(H_{J^r}) = 
\sum_{J \succcurlyeq I}(-1)^{\ell(J) - \ell(I)} \sum_{K \preccurlyeq J^{r r}} (-1)^{\ell(K)} H_K = \\
& = \sum_{K} (-1)^{\ell(K)}  H_{K}  
\sum_{\substack{J\succcurlyeq K \\ J\succcurlyeq I} }(-1)^{\ell(J)-\ell(I)} =
(-1)^{|I|}  
\sum_{K\succcurlyeq I^c }(-1)^{\ell(K)-\ell(I^c)} H_{K} =  (-1)^{|I|} R_{I^c} .
\end{align*}
Above we sum up with respect to all sentences $K$ having the same maximal word $w(I)$, and the inner sum is over sentences $J$ which coarsen simultaneously $K$ and $I$. This sum doesn't vanish only when this poset of sentences is trivial what happens only if  $K$ coarsens the complement of $I$. In this case $J$ is the single-word sentence, thus $\ell(J)=1$. To conclude the calculation we notice that  $\ell(I) + \ell(I^c) = |I|-1$. 
\end{proof}

\section{Formal series of trees and coloured non-commutative power sum functions} \label{sec:restricted}
Up to now we considered the duality problem for infinite-dimensional Hopf algebras in the graded case only. 
Another option to tackle the problem is to define~\cite{Abe,Sweedler} the \emph{restricted (or Sweedler's) dual} of $\mathcal{H}$ which is the subspace $\mathcal{H}^\circ \subset \mathcal{H}^*$ consisting of all linear maps $f$ that satisfy one of equivalent conditions:
\begin{enumerate}
	\item $\Delta_{\mathcal{H}^*} (f) \in \mathcal{H}^* \otimes \mathcal{H}^*$,
	\item $\mathrm{ker} (f) $ contains an ideal (left, right or two-sided) of $\mathcal{H}$ that has finite codimension.
\end{enumerate}

Define a left action $\rightharpoonup$ of $\mathcal{H}$ on $\mathcal{H}^*$ as the transpose of right multiplication on $\mathcal{H}$
\begin{equation*}
\langle a \rightharpoonup f, b \rangle = \langle f, b a \rangle \qquad f\in \mathcal{H}^*, \quad a,b \in \mathcal{H}.
\end{equation*}
Then $\mathcal{H} \rightharpoonup f$ is a subspace of $\mathcal{H}^*$, and the condition $f\in \mathcal{H}^\circ$ is equivalent to
\begin{enumerate}
	\item[(3)] $\dim ( \mathcal{H} \rightharpoonup f) < \infty$.
\end{enumerate}
\begin{Rem}
	One can define also a right action $\leftharpoonup$ of $\mathcal{H}$ on $\mathcal{H}^*$ as the transpose of left multiplication on $\mathcal{H}$. Then condition (3) can be equivalently stated as finite-dimensionality of $f \leftharpoonup \mathcal{H}$ or  finite-dimensionality of  $ \mathcal{H} \rightharpoonup f \leftharpoonup \mathcal{H}$.
\end{Rem}

\subsection{Formal series of ROC trees}
In this Section we consider power series of trees as the linear dual to space of rooted ordered coloured (by $A$) trees $\Bbbk T_A$. 
A \emph{formal tree series} $F$ is function $T_A \to \Bbbk$ extended to $\Bbbk T_A $ by linearity.
The image by $F\in (\Bbbk T_A)^*$ of a tree $ t\in T_A$ is denoted by $\langle F,t \rangle $ and is called the coefficient of $t$ in $T$. The support of $F$ is the subset of $T_A$
\begin{equation}
\mathrm{supp}(F) = \{ t\in T_A | \langle F,t \rangle \neq 0 \}.
\end{equation}
Polynomials $\Bbbk T_A   \subset (\Bbbk T_A)^*$ are embedded naturally as series with finite support.  Usually one writes 
\begin{equation}
F = \sum_{t\in T_A} \langle F, t\rangle t, 
\end{equation}
remembering that the sum
\begin{equation}
\Bbbk T_A  \ni P \mapsto \langle F,P \rangle =  \sum_{ t \in T_A} \langle F,t \rangle \langle P, t \rangle ,
\end{equation}
has a finite support. 
\begin{Rem}
	Notice that, by the standard coding of trees by (coloured) Dyck words~\cite{Stanley}, any such series of trees can be interpreted as a series of words within theory of non-commutative power series~\cite{Salomaa}. As it was mentioned in Corollary~\ref{cor:NSym-CF}, our description of the algebra $\mathrm{NSym}_A$ can be stated in terms of a certain context-free language.
\end{Rem}
\begin{Rem}
	If $\Bbbk$ is equipped with discrete topology, then the set of formal tree series 
can be equipped with the product topology. A sequence of its elements converges only if for each tree the corresponding coefficient stabilizes.
\end{Rem}
Actually, two products of such series are well defined:
\begin{itemize}
	\item the extension of the concatenation product $" . "$ of trees (i.e. the Cauchy product of series)
\begin{equation}
F . G = \sum_{t\in T_A} \left( \sum_{t^\prime . t^{\prime\prime}} \langle F, t^\prime \rangle \langle G,  t^{\prime\prime} \rangle  \right) t ;
\end{equation}
	\item the extension of the grafting product $\Tshuffle$ of trees
\begin{equation}
F \Tshuffle G =  \sum_{t^\prime , t^{\prime\prime} \in T_A} \langle F, t^\prime \rangle \langle G,  t^{\prime\prime} \rangle \;   t^\prime \Tshuffle t^{\prime\prime} .
\end{equation}
\end{itemize}  

The tensor product $F\otimes G \in (\Bbbk T_A)^* \otimes (\Bbbk T_A)^*$ of two series reads
\begin{equation}
F\otimes G = \sum_{t,s\in T_A} \langle F, t \rangle \langle G ,s \rangle  t\otimes s .
\end{equation}
The deconcatenation coproduct extended from tree polynomials to series 
\begin{equation}
\delta(F) =  \sum_{t,s\in T_A} \langle F,t.s\rangle t\otimes s ,
\end{equation}
is in general an element of $(\Bbbk T_A \otimes \Bbbk T_A )^* $.
A series $F$ which allows for finite decomposition
\begin{equation} \label{eq:delta-S} 
\delta(F) =  \sum_{i=1}^r G_i \otimes H_i \in (\Bbbk T_A)^* \otimes (\Bbbk T_A)^*,
\qquad \text{with} \qquad G_i, H_i \in   (\Bbbk T_A )^*,
\end{equation}
is an element of the restricted dual $  (\Bbbk T_A, \: .\: , \Delta)^\circ$. The pruning coproduct of a tree series can be defined analogously
\begin{equation}
\Delta(F) =  \sum_{t,s\in T_A} \langle F,t\Tshuffle s\rangle t\otimes s .
\end{equation}

Let us describe a distinguished example of two of such tree series. By $F^*$ let us denote the characteristic series of the set of all ROC trees, 
i.e. $F^*=\sum_{t \in T_A} t$, and by $F=\sum_{t \in T_A^\prime} t$ denote the characteristic series of the subset $T^\prime_A$ of planted trees. Relations between these series can be written down as follows:

(i) an arbitrary planted tree is obtained by action of the operator $B_i^+$, $i=1,2,\dots , |A|$, on the corresponding tree
	\begin{equation} \label{eq:fT'}
	F = \sum_{i=1}^{|A|} B_i^+(F^*) \; ,
	\end{equation}
	
(ii)  any non-trivial tree can be uniquely decomposed into concatenation product of planted trees (this justifies our notation)
	\begin{equation} \label{eq:fT-fT'-s}
	F^* = \ttree + F + F  . F + F . F  . F + \dots = \ttree +  F. F^* \; .
	\end{equation}
\begin{Rem}
	Equations \eqref{eq:fT'}-\eqref{eq:fT-fT'-s} expresses the standard grammar rules of the $|A|$-th Dyck language~\cite{Salomaa}.
\end{Rem}
\begin{Rem}
	By combining equations~\eqref{eq:fT'} and \eqref{eq:fT-fT'-s} we obtain a single equation for series $F$ 
\begin{equation}
F = \sum_{i=1}^{|A|} B_i^+ \left(\ttree + \sum_{k=1}^\infty F^{k} \right) 
\end{equation}	
in  the form of the combinatorial Dyson--Schwinger equation~\cite{Broadhurst-Kreimer,Foissy-DS}.
\end{Rem}
\begin{Prop} \label{prop:delta-F-F'}
	The series $F$ and $F^*$ are elements of the restricted dual $(\Bbbk T_A, .\,, \Delta)^\circ$ of the Hopf algebra of rooted ordered coloured trees, in particular $F$ is primitive element
	\begin{equation} \label{eq:delta-ft'}
	\delta(F)  = \ttree \otimes F + F \otimes \ttree, 
\end{equation}
and $F^*$ is group-like
\begin{equation} \label{eq:delta-ft}
	\delta(F^*)  = F^* \otimes F^* \; .
	\end{equation}
\end{Prop}
\begin{proof}
	Equation \eqref{eq:delta-ft'} follows from the analogous result valid for any planted tree $t^\prime \in T_A^\prime$
	\begin{equation*} 
	\delta(t^\prime) = \ttree \otimes t^\prime + t^\prime \otimes \ttree .
	\end{equation*}
	When $t^\prime \in T_A^\prime$ is a planted tree and $t\in T_A$ is an arbitrary tree then by Corollary~\ref{cor:delta-.}
	\begin{equation} \label{eq:delta-t'.t}
\delta(t^\prime . t) = \ttree \otimes (t^\prime . t) + (t^\prime \otimes \ttree) . \delta (t) ,
	\end{equation}
	which by linearity leads to the following equation on the level of the corresponding series 
	\begin{equation}
	\delta(F . F^* ) = \ttree \otimes (F . F^*) + (F \otimes \ttree) . \delta(F^*).
	\end{equation}
	From equation \eqref{eq:fT-fT'-s} we obtain the relation
	\begin{equation}
	\delta(F^*) 	= \ttree \otimes F^* + (F  \otimes \ttree) . \delta(F^*),
	\end{equation}
	which, when solved for $\delta(F^*)$, gives  \eqref{eq:delta-ft}.
\end{proof}
\begin{Cor}
	Any series, whose support is a subset of planted trees $T^\prime_A$ is primitive element of the restricted dual.
\end{Cor}
\begin{Rem}
	Notice that Corollary~\ref{cor:delta-.} implies the following matching condition between the cut comultiplication~$\delta$ and the Cauchy product of two series $G$ and $H$ of trees in the restricted dual
\begin{equation}
\delta(G.H) = \delta(G) . (\ttree \otimes H) + (G\otimes \ttree).\delta(H) - G\otimes H.
\end{equation}	
\end{Rem}

\subsection{Formal series of tall trees}
By $H^*=\sum_I I$ denote the characteristic series of the set of coloured tall trees (indexed by sentences over $A$) and by $H = \sum_{w\in A^*\setminus\{1\}} (w)$ denote the characteristic series  of planted coloured tall trees (indexed by single word sentences). The planting operator $B_a^+$, $a\in A$, acts on single word sentence $(w)$ by forming the single word sentence $(wa)$. Properties of the series can be stated as follows:

(i) an arbitrary planted tall tree is obtained by action of the operator $B_i^+$, $i=1,2,\dots ,|A|$ on the trivial tree or on the corresponding smaller planted tall tree
\begin{equation}
H = \sum_{i=1}^{|A|} B_i^+(\ttree + H) \label{eq:fH'} \; ,
\end{equation}
(ii) any non-trivial tall tree can be uniquely decomposed into concatenation of planted tall trees
\begin{equation}
H^* = \ttree + H + H  . H  + H  . H  . H + \dots = \ttree + H . H^* \; . \label{eq:fH-fH'}
\end{equation}
\begin{Rem}
	Equation \eqref{eq:fH'} is of the combinatorial Dyson--Schwinger form. 
\end{Rem}
	
The following result can be proven in the same way like the previous Proposition~\ref{prop:delta-F-F'}.
\begin{Prop}
With respect to the deconcatenation coproduct $\delta$ the series $H$ is primitive element of the restricted dual
\begin{equation} \label{eq:delta-fH'}
\delta(H) = \ttree \otimes H + H \otimes \ttree, 
\end{equation}
and the series $H^*$ is group-like
\begin{equation} \label{eq:delta-fH}
\delta(H^*)  = H^* \otimes H^* \; .
\end{equation}	
\end{Prop}
\begin{Rem}
	Actually, for action of the coproduct $\delta$ we should write the series $H$ and $H^*$ in terms of the coloured monomial quasi-symmetric functions
	\begin{gather}
	H_Q = \sum_{w\in A^* \setminus \{ 1\}} M_{(w)}, \qquad 
	H_Q^* = \sum_{I} M_I,\\
	\delta(H_Q) = 1 \otimes H_Q + H_Q \otimes 1, \qquad
	\delta(H_Q^*)  = H_Q^* \otimes H_Q^* \; .
	\end{gather} 
	i.e. as elements of the restricted dual $(\mathrm{NSym}_A)^\circ$. 
\end{Rem}
Finally let us present formulas for antipodes of the above series. Because $H_Q$ is primitive element then 
\begin{equation}
S^*(H_Q) = -H_Q,
\end{equation}
but the antipode of the sum of all coloured monomial quasi-symmetric functions takes also particularly simple form.
\begin{Cor}
The action of the antipode $S^*$ on the series $H_Q^*$ reads
\begin{equation}
S^*(H_Q^*) = \sum_{w\in A^*} (-1)^{|w|} M_{(w^c)},
\end{equation}
i.e. is the signed sum of monomial functions indexed by minimal compositions (i.e. by sentences built from single-letter words).
\end{Cor}
\begin{proof}
	Decompose the series $H_Q^*$ into parts indexed by compositions with the same maximal word 
\begin{equation}
H_Q^* = \sum_{I}M_I = \sum_{w\in A^*} F_{(w)},
\end{equation}
where $F_{(w)}$ is the fundamental function indexed by the corresponding single-word composition,
and apply Proposition~\ref{prop:S*F}.
\end{proof}

Because bases of both algebras $\mathrm{NSym}_A$ and $\mathrm{QSym}_A$ are indexed by sentences one can apply also the pruning coproduct $\Delta$ on series of tall trees. In particular, we will show that the series $H$ can be considered as an element of the restricted dual $(\mathrm{QSym}_A)^\circ$. From now on we use the notation of the theory of the non-commutative coloured symmetric functions.
\begin{Prop}
	The pruning coproduct of the series $ H = \sum_{w\in A^+} H_{(w)}$, where $A^+ = A^*\setminus \{1\}$
reads
\begin{equation}
\Delta(H) = 1 \otimes H +  H \otimes 1 + 
H \otimes H.
\end{equation}
\end{Prop}
\begin{proof}
By Corollary~\ref{cor:Delta-w}, for arbitrary one word sentence $(w)$ we have
\begin{equation}
\Delta(H_{(w)}) = \sum_{uv=w}H_{(u)}\otimes H_{(v)}.
\end{equation}
Therefore, summing up with respect to the arbitrary prefix $(u)$ first, we can write
\begin{equation*}
\Delta \left( \sum_{w\in A^+} H_{(w)} \right) = \sum_{w\in A^+} \sum_{uv=w}H_{(u)}\otimes H_{(v)} =
\sum_{u \in A^*} H_{(u)}\otimes \left( \sum_{v\in u^{-1}A^+} H_{(v)} \right),
\end{equation*}
where for any subset $L\subset A^*$ by definition~\cite{Sakarovitch}
\begin{equation*}
u^{-1}L =  \{ v \, | \,  uv \in L \}.
\end{equation*} 
Because
\begin{equation*}
u^{-1}A^+  = \begin{cases}
A^+ & u=1, \\ A^* & u\in A^+ ,
\end{cases} 
\end{equation*}
then we obtain
\begin{equation*}
\Delta(H) = 1 \otimes H + H \otimes  (1 + H).
\end{equation*}
\end{proof}
\begin{Cor}	\label{cor:H-gl}
The series $1 + H = \sum_{w\in A^*} H_{(w)}$ is group like element of the restricted dual $(\mathrm{QSym}_A)^\circ$.
\end{Cor}

\subsection{Non-commutative coloured power sum symmetric functions} \label{sec:power-sum}
Define the power series 
\begin{equation}
P = \log (1 + H) = \sum_{n=1}^\infty \frac{(-1)^{n-1}}{n} H^n,
\end{equation}
which exists because $H$ has vanishing constant term~\cite{Salomaa}.
\begin{Prop}
In the basis of complete homogeneous functions indexed by compositions we have
\begin{equation} \label{eq:P-H}
P= \sum_{I \neq 1} \frac{(-1)^{\ell(I)-1}}{\ell(I)} H_I, 
\end{equation}
moreover the series $P$ is primitive with respect to the pruning coproduct
\begin{equation} \label{eq:Delta-P}
\Delta(P) = 1\otimes P + P\otimes 1 .
\end{equation}
\end{Prop}
\begin{proof}
	The first part follows directly from the definition of $P$, where we recall that the length $\ell(I)$ of the sentence $I$ is the number of its words.	
	For the second part we provide two proofs. The first one follows the corresponding reasoning~\cite{MalvenutoReutenauer} in the monochromatic case. The second proof is of elementary combinatorial nature. 

I. By linearity and morphism property of $\Delta$ with respect to the concatenation product, and using Corollary~\ref{cor:H-gl} we have	
	\begin{multline*}
	\Delta(P) =  \log \left[ \Delta (1 + H ) \right] = \log \left[(1+H)\otimes  (1+H)\right] = \log \left[ \left( (1+H)\otimes  1 \right) . \left( 1\otimes  (1+H)\right) \right] = \\
	=  \log \left[ (1+H)  \otimes  1 \right] + \log \left[ 1\otimes  (1+H)\right] = 
	\log \left[ 1+H  \right] \otimes  1 + 1\otimes    \log \left[1+H\right] = P\otimes 1 + 1 \otimes P, \; \; \; 
	\end{multline*}
where we also used the standard property of logarithm for commuting factors.

II. By the basic coproduct formula \eqref{eq:Delta-HI} applied to equation~\eqref{eq:P-H} we can see that 
\begin{equation}
\Delta(P) = P\otimes 1 + 1 \otimes P + \sum_{J,K\neq 1} c_{JK} \, H_J \otimes H_K,
\end{equation}
where the coefficient $c_{JK}$ equals
\begin{equation}
c_{JK} = \sum_I \frac{(-1)^{\ell(I)-1}}{\ell(I)} ,
\end{equation}
where we sum with respect to the sentences $I$ which give $J\otimes K$ upon action of the coproduct, i.e.~$I$~is a~summand in the quasi-shuffle $J\Qshuffle K$. Interpretation of such terms as special paths in the lattice $\ell(J) \times \ell(K)$, and application of Proposition~\ref{prop:q-shuffle-paths-i} implies that the weighted alternating sum we are looking for equals (without losing generality we assume $\ell(K)\geq \ell(J)$)
\begin{equation*}
\frac{(-1)^{\ell(J) + \ell(K)-1}}{\ell(J)}\sum_{i=0}^{\ell(J)} (-1)^i \binom{\ell(J) }{i} \binom{\ell(K)+\ell(J)-i-1}{\ell(J)-1},
\end{equation*}
which vanishes by standard application of the inclusion-exclusion principle.
\end{proof}

Finally, we define the coloured non-commutative analogs of the power sum symmetric functions. By splitting series $P$ into parts with the same maximal words, see Corollary~\ref{cor:maximal-word}, we obtain from equation~\eqref{eq:P-H}
\begin{equation}
P = \sum_{w\in A^+} P_{(w)}, \qquad P_{(w)}=\sum_{I\preccurlyeq (w)}\frac{(-1)^{\ell(I)-1}}{\ell(I)} H_I. 
\end{equation}
Here $|w|P_{(w)}$ are coloured analogs of the non-commutative power sums of the second kind defined in~\cite{Gelfand-NSym}. For trivial sentence define $P_1 = 1$, and for any non-empty sentence $I=(w_1,w_2,\dots , w_{\ell(I)})$ define 
\begin{equation}
P_I = P_{(w_1)} P_{(w_2)}  \dots P_{(w_{\ell(I)})}.
\end{equation}
\begin{Prop}
	Functions $P_I$ indexed by sentences form a linear basis of $\mathrm{NSym}_A$, in particular
\begin{equation} \label{eq:H-P}
H_{(w)} = \sum_{I\preccurlyeq (w)}\frac{1}{\ell(I)!} P_I .
\end{equation}
\end{Prop}
\begin{proof}
	By the standard relation between exponential and logarithm, valid also for formal non-commuting series, we have
	\begin{equation}
	1+\sum_{w\in A^+} H_{(w)} = \exp(P) = \sum_{n=0}^\infty \frac{1}{n!} P^n =
	\sum_{n=0}^\infty \frac{1}{n!} \left( \sum_{w\in A^+} P_{(w)} \right)^n = 	
	1 + \sum_{I\neq 1} \frac{1}{\ell(I)!} P_I,
	\end{equation}
	and formula \eqref{eq:H-P} follows from splitting of both sides into sentences with the same maximal word.
\end{proof}
\begin{Ex}
	For $w=ab$ we have
\begin{equation*}
P_{(ab)} = H_{(ab)} - \frac{1}{2} H_{(a)}H_{(b)}, \qquad \Delta(P_{(ab)}) = 1 \otimes P_{(ab)} + P_{(ab)} \otimes 1 + \frac{1}{2} \left( H_{(a)} \otimes H_{(b)} - 
H_{(b)} \otimes H_{(a)} \right).
\end{equation*}
\end{Ex}
As the above example demonstrates, contrary to the monochromatic/unary case the coloured power sum functions are in general not primitive elements of the Hopf algebra $\mathrm{NSym}_A$.
However, by splitting equation \eqref{eq:Delta-P} into homogeneous parts we obtain the following weaker result, which provides infinite number of primitive elements of the algebra.
\begin{Cor}
	For $n\in\mathbb{N}$ define $P_n = \sum_{ |w|=n} P_{(w)}$ then
\begin{equation} 
\Delta(P_n) = 1\otimes P_n + P_n \otimes 1 .
\end{equation}	
\end{Cor}
\section{Conclusion} 
We defined new generalization $\mathrm{NSym}_A$ and $\mathrm{QSym}_A$ of the Hopf algebras of non-commutative symmetric and quasi-symmetric functions. In our extension bases in both algebras are indexed by sentences over finite alphabet $A$ (the set of colours), and our results reduce to the classical ones for unary alphabet $|A|=1$ (the monochromatic reduction). We presented corresponding analogs of the most pertinent structural elements of the original theory including description of bases of complete homogeneous, elementary, monomial, fundamental, ribbon Schur and power sum functions.
It is interesting that, contrary to the monochromatic/unary case, both algebras are non-commutative and non-cocommutative. We found also realization of the algebra $\mathrm{QSym}_A$ in terms of power series of bounded degree in partially commuting variables, which justifies its name as coloured quasi-symmetric functions. 

In our approach the algebra $\mathrm{NSym}_A$ is described as Hopf subalgebra of rooted ordered coloured trees. We study also formal series of such trees within the setting of the restricted duals. This new aspect of the theory deserves deeper studies in relation to renormalization procedure in quantum field theory, non-commutative integrable systems and context-free languages, and will be developed in another publication. 
In the literature there are known several generalizations of the non-commutative symmetric and quasi-symmetric functions. We strongly believe that the generalization proposed in our paper, being natural and structurally close to the original theory, will be useful in studying problems in combinatorics and physics. 

\bibliographystyle{amsplain}

\end{document}